\documentclass[11pt]{article}
\usepackage{amsmath,amssymb,eucal}
\usepackage[dvips]{graphicx}

\newtheorem{theorem}{\bf Theorem}[section]

\newtheorem{lemma}[theorem]{\bf Lemma}

\newtheorem{definition}[theorem]{\bf Definition}
\newtheorem{notation}[theorem]{\bf Notation}

\newtheorem{question}[theorem]{\bf Question}

\numberwithin{equation}{section}

\def\qed{$\Box$}

\begin{document}
\title{
Divide knot presentation of sporadic knots of \\
Berge's lens space surgery 
}
\author{Yuichi YAMADA}
%\date{\today}
\date{October 10, 2012}
\footnotetext[0]{%
}
\maketitle
%%%%%%%%%%%%%%%%
\begin{abstract}{
Divide knots and links, defined by A'Campo in
the singularity theory of complex curves, is a method to 
present knots or links by real plane curves.
The present paper is a continuation of the author's previous result that 
every knot in the major subfamilies of Berge's lens space surgery (i.e., 
knots yielding a lens space by Dehn surgery) is presented by an L-shaped 
curve as a divide knot.
In the present paper, L-shaped curves are generalized and 
it is shown that every knot in the minor subfamilies, called sporadic examples, 
of Berge's lens space surgery is presented by a generalized L-shaped curve as a divide knot. 
A formula on the surgery coefficients and the presentation is also generalized. 
}
\end{abstract}
%%%%%%%%%%%%%%%%

%%%%%%%%%%%%%%%%%%%%%%%%%%%%%%%%%%%%%%%%%%
%%%                   Section 1
%%%%%%%%%%%%%%%%%%%%%%%%%%%%%%%%%%%%%%%%%%
\section{Introduction}~\label{sec:intro}
If $r$ Dehn surgery on a knot $K$ in $S^3$ yields the
lens space $L(p,q)$, we call the pair $(K, r)$ a {\it lens space surgery},
and we also say that $K$ admits a lens space surgery,
and that $r$ is the {\it coefficient} of the lens space surgery.
The task of classifying lens space surgeries, especially 
knots that admit lens space surgeries has been
a focal point in low-dimensional topology.
In 1990, Berge \cite{Bg}
pointed out a \lq\lq mechanism\rq\rq \ %
of known lens space surgery,
that is, {\it doubly-primitive knots} in the Heegaard surface of genus $2$.
Berge also gave a conjecturally complete list of such knots,
described them by Osborne--Stevens's
\lq\lq R-R diagrams\rq\rq \ %
in \cite{OS}, and classified such knots into three families,
and into 12 types in detail:
\begin{enumerate}
%%%
\item[(1)] {\it Knots in a solid torus (Berge--Gabai knots)} : TypeI, II, ... and VI 
(Berge \cite{Bg2}) 
\par
Dehn surgery along a knot in a solid torus whose resulting
manifold is also a solid torus. TypeI consists of torus knots.
TypeII consists of $2$-cable of torus knots.
%%%
\item[(2)] {\it Knots in genus-one fiber surface} : TypeVII and VIII
(see Baker \cite{Ba, Ba3} and also \cite{Y1, Y5})
%%%
\item[(3)] {\it Sporadic examples (a), (b), (c) and (d)} : 
TypeIX, X, XI and XII, respectively.
\end{enumerate}
Their surgery coefficients are also decided.
They are called {\it Berge's lens space surgeries}.
The numbering VII--XII (after VI) are used 
by Baker in \cite{Ba2, Ba3}.
It is conjectured by Gordon \cite{go1, go2} that
every knot of lens space surgery is a doubly-primitive knot.

In the present paper, we are concerned with the minor family (3).
It is known
that TypeIX and TypeXII (Berge's (a) and (d)) are related, 
and 
that TypeX and TypeXI (Berge's (b) and (c)) are related.
Thus our targets are TypeIX and TypeX.
%
%
%%%%%%%%%%%%%%%%%%%%%%%%%%%%%%%%%%%%
\begin{notation}
Throughout the paper, we let Type${\mathcal X}$
denote either TypeIX or TypeX, i.e., ${\mathcal X} =$ IX or X.
Knots in Type $\mathcal{X}$ are parametrized
by an integer $j$ with $j \not= 0,-1$ (\cite{Bg}),
thus we call the knots $k_{IX}(j)$ for TypeIX,
and $k_{X}(j)$ for TypeX, respectively.
For the precise construction of the knots $k_{\mathcal{X}}(j)$, 
see Section~\ref{sec:sporadic}.
\end{notation}
%%%%%%%%%%%%%%%%%%%%%%%%%%%%%%%%%%%%
%
%
Berge's original classification (a)--(d) of sporadic examples in \cite{Bg} 
was 
\begin{center}
\begin{tabular}{ll}
(a) \ $k_{IX}(j)$ with $j > 0$, 
&
(b) \ $k_{X}(j)$ with $j > 0$, \\ 
(c) \ $k_{X}(j)!$ with $j < -1$, 
&
(d) \ $k_{IX}(j)!$ with $j < -1$,
\end{tabular}
\end{center}
respectively, see Deruelle--Miyaszaki--Motegi's recent works \cite{DMM, DMM2}. 
Table~\ref{tbl:sporad} is a list of some data of 
the knots $k_{\mathcal{X}}(j)$:
the coefficients $p$ of their lens space surgeries of $k_{\mathcal{X}}(j)$,
the second parameter $q$ of the resulting lens space $L(p,q)$
and the genus of $k_{\mathcal{X}}(j)$, which depends on the sign of $j$.
Our convention about orientations of lens spaces is 
\lq\lq $p/q$ Dehn surgery along an unknot is $-L(p,q)$\rq\rq.
%
%
%%%%%%%%%%%%%%%%%
\begin{table}[h]
\begin{center}
\begin{tabular}{|l|c|c|c|c|}
\hline
knot  & $p$ ($=$ coefficient) & $q$ of $L(p,q)$ & genus ($j>0$)& genus ($j< -1$)\\
\hline
$k_{IX}(j)$ & $22 j^2 + 9j +1$ & $-(11j+2)^2$ & $11j^2$ & $11 j^2 + 9j +2$ \\
\hline
$k_{X}(j)$ & $22 j^2 + 13j +2$ & $-(11j+3)^2$ & $11j^2+2j$ & $11j^2 + 11j +3$\\
\hline
\end{tabular}
\caption{Data on Sporadic knots}
\label{tbl:sporad}
\end{center}
\end{table}
%%%%%%%%%%%%%%%%%
%
%

\medskip

%
%
%%%%%%%%%%%%%%%%%
\begin{figure}[h]
\begin{center}
\includegraphics[scale=0.4]{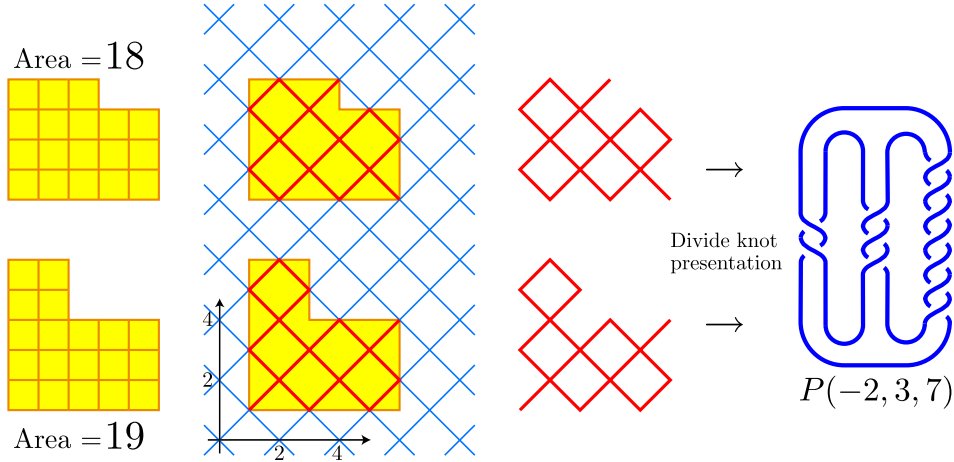}\\
\caption{Divide knot presentation of $P(-2,3,7)$ (see \cite{Y4, Y1})}
\label{fig:IntroC}
\end{center}
\end{figure}
%%%%%%%%%%%%%%%%%
%
%
The theory of A'Campo's {\it divide knots and links}
came from singularity theory of complex curves.
A {\it divide} is originally
a relative, generic immersion of a 1-manifold
in the unit disk in $\mathbb{R}^2$, see Section~\ref{sec:curve}.
A'Campo \cite{A1, A2, A3, A4}
formulated the way to associate to each divide $P$ a link
$L(P)$ in $S^3$.
In the present paper, we regard a PL (piecewise linear) plane curve as a divide
by smoothing the corners and controlling the size. 
Let $X$ be the $\pi /4$-lattice
defined by $\{ (x,y) \vert \cos \pi x = \cos \pi y \}$
in $xy$-plane. 
In this paper, we are interested in plane curves constructed 
as intersection of $X$ and a region.
See Figure~\ref{fig:IntroC}, which was the starting examples
of the author's project.
Two L-shaped curves of the form $X\cap \mathcal{L}$
present a same knot, the pretzel knot of type $(-2,3,7)$, 
as divide knots.
Its $18$-surgery and $19$-surgery are lens spaces 
(by Fintushel--Stern \cite{FS}).
Note that the areas of $\mathcal{L}$ are equal to 
the coefficients of the lens space surgeries.
They have different mechanism of lens space surgeries:
$18$-surgery is in TypeIII, $19$-surgery is in TypeVII. 

\medskip
Our question is
%
%
%%%%%%%%%%%%%%%%%%%%%%%%%%%%%%%%%%%%
\begin{question}
Is every knot of (Berge's) lens space surgeries
a divide knot? 
\end{question}
%%%%%%%%%%%%%%%%%%%%%%%%%%%%%%%%%%%%
%
%
%
%
%%%%%%%%%%%%%%%%%
\begin{figure}[h]
\begin{center}
\includegraphics[scale=0.4]{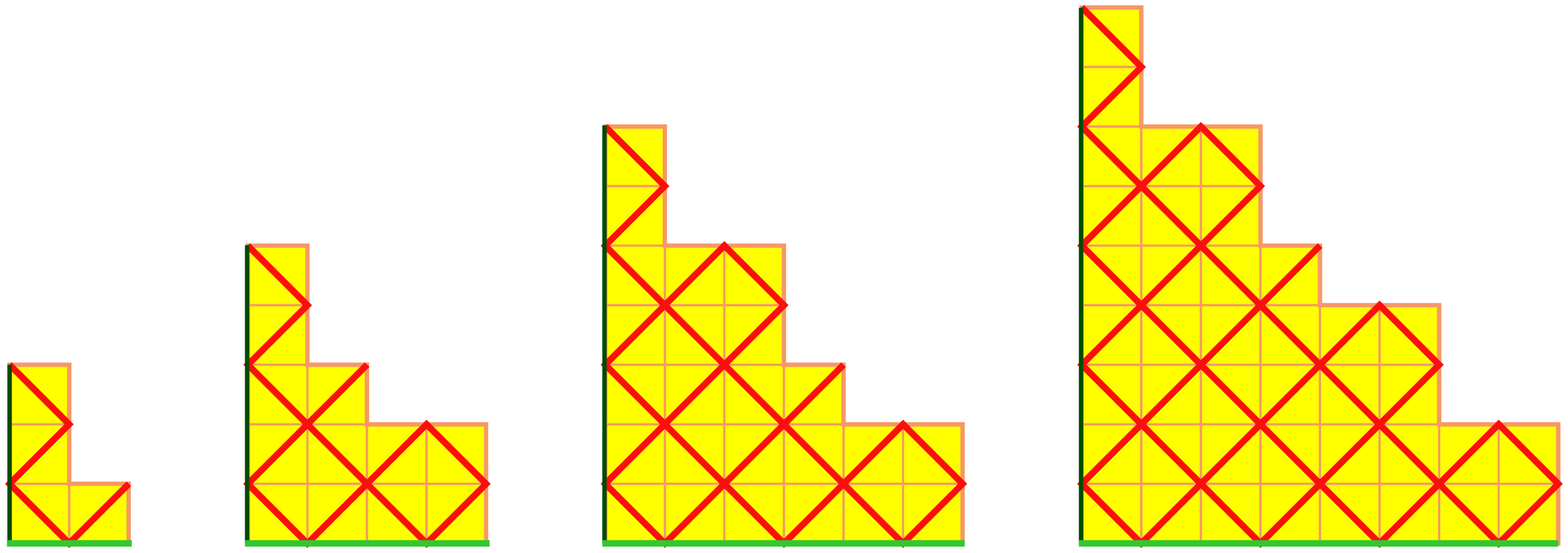}\\
Case $j > 0$ : \ $P(1), P(2), P(3), P(4)$ \\
\bigskip
%
%\caption{$P(j)$ with $j>0$ (ex. $P(1), P(2), P(3), P(4)$)}
\includegraphics[scale=0.4]{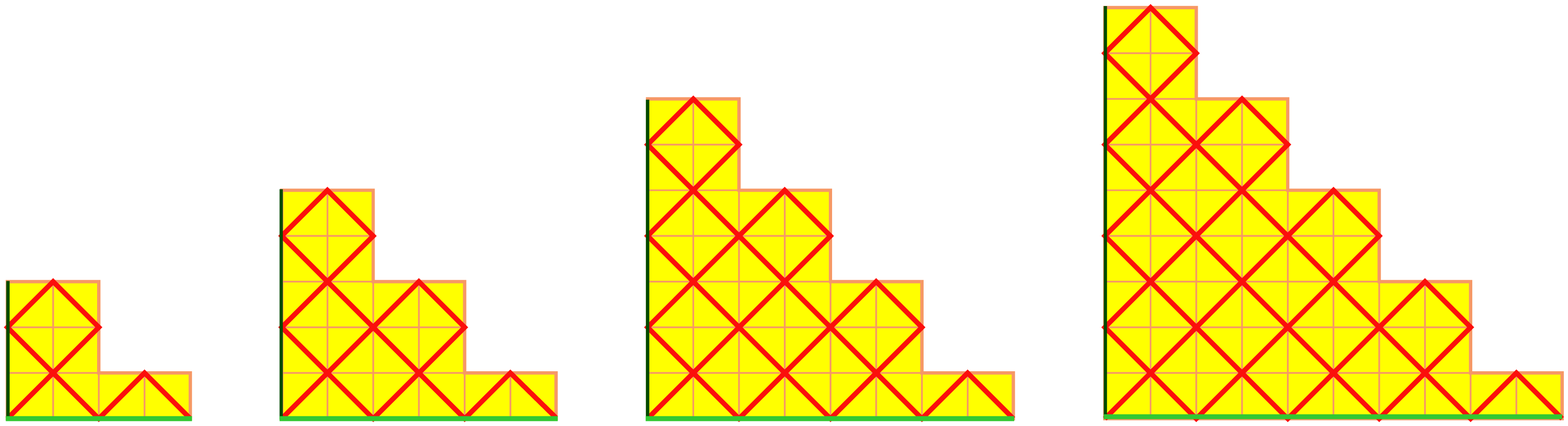}\\
Case $j < -1$ : \ $P(-2), P(-3), P(-4), P(-5)$
\caption{Plane curves $P(j)$, see Definition~\ref{def:p(j)}
for the precise definition}
\label{fig:SporadPN}
\end{center}
\end{figure}
%%%%%%%%%%%%%%%%%
%
%
The purpose of the present paper is to give
divide knot presentations of knots of sporadic examples.
We show some plane curves to state the main results,
see $P(j)$ in Figure~\ref{fig:SporadPN}.
Each curve $P(j)$ is constructed 
as intersection of $X$ and a region that consists of 
some rectangles sharing the left bottom corner.
We call such a curve {\it a generalized L-shaped curve}.
The precise definition will be given in Section~\ref{sec:curve} and Section~\ref{sec:sporadic}.

Next, we define plane curves $P_{\mathcal{X}}(j)$
from $P(j)$ by {\it adding a square}
twice in the sense \cite[Lemma~4.2]{Y4}.
%
%
%%%%%%%%%%%%%%%%%
\begin{figure}[h]
\begin{center}
\includegraphics[scale=0.4]{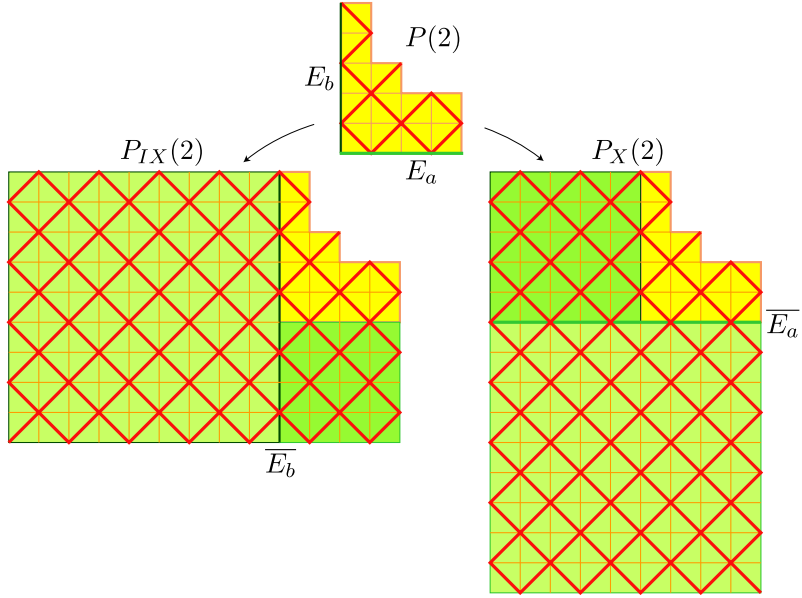}
\caption{Add squares twice to get $P_{\mathcal{X}}(j)$ \ 
(ex. $P_{IX}(2)$ and $P_{X}(2)$)}
\label{fig:SporadEx}
\end{center}
\end{figure}
%%%%%%%%%%%%%%%%%
%
%
%
%
%%%%%%%%%%%%%%%%%
\begin{definition}[Plane curves $P_{\mathcal{X}}(j)$]~\label{def:PXj}
We let $E_a$ (and $E_b$, respectively) denote the bottom edge 
(and the left edge) of the region of $P(j)$, see Subsection~3.1 for the precise definition.
For an integer $j$ with $j \not=0,-1$, depending on whether $\mathcal{X} =$ IX or X,
we construct a region of $P_{\mathcal{X}}(j)$ as follows, 
see Figure~\ref{fig:SporadEx}.
\begin{enumerate}
\item[]
[TypeIX] \
We add a square along the bottom edge $E_a$ first,
and add another square along the lengthened left edge
$\overline{E_b}$. 
\item[] 
[TypeX] \
We add a square along the left edge $E_b$ first,
and add another square along the lengthened bottom edge
$\overline{E_a}$. 
\end{enumerate}
\end{definition}
%%%%%%%%%%%%%%%%%
%
%
We remark that, by the first square addition
along an edge $E_a$ (or $E_b$, respectively),
the other edge $E_b$ (or $E_a$) is 
lengthened as $\overline{E_b}$ (or $\overline{E_a}$).
The second square is added along the lengthened one.
By $l(E)$ we denote the length of the edge $E$. Then,
\[
l(E_a) = \vert 2j \vert, \ 
l(E_b) = \vert 2j+1 \vert
\text{ and }
l(\overline{E_b}) = 
l(\overline{E_b}) = l(E_a) + l(E_b).
\]

\medskip

Our main theorem is 
%
%
%%%%%%%%%%%%%%%%%%%%%%%%%%%%%%%%%%%%
\begin{theorem}~\label{thm:main}
Up to mirror image,
every knot in TypeIX and TypeX in 
Berge's list of lens space surgery is a divide knot. 
More precisely, 
the plane curve $P_{\mathcal{X}}(j)$
presents the knot $k_{\mathcal{X}}(j)$ as a divide knot.
\end{theorem}
%%%%%%%%%%%%%%%%%%%%%%%%%%%%%%%%%%%%
%
%
We will also show
%
%%%%%%%%%%%%%%%%%%%%%%%%%%%%%%%%%%%%
\begin{lemma}~\label{lem:Pj}
The plane curve $P(j)$
presents the torus knot $T(j, 2j+1)$ as a divide knot.
\end{lemma}
%%%%%%%%%%%%%%%%%%%%%%%%%%%%%%%%%%%%
%
%
%%%%%%%%%%%%%%%%%%%%%%%%%%%%%%%%%%%%
\begin{lemma}~\label{lem:Pm}
{\rm (see \cite{DMM} on the knots)}
\ We let $P_m(j)$ denote a plane curve obtained 
by a square addition along $E_a$ from $P(j)$, which appears
in the process [TypeIX] to construct $P_{IX}(j)$
in Definition~\ref{def:PXj}:
\[
P(j) 
\ \rightarrow \
P_m(j)
\ \rightarrow \
P_{IX}(j).
\]
The plane curve $P_m(j)$ presents the cable knot 
$C(T(2,3); j, 6j + 1)$,
regarded as $C(T(2,3); \vert j \vert, 6\vert j\vert - 1)$ if $j<-1$,
of the trefoil as a divide knot:
\[
T(j,2j+1) 
\ \rightarrow \
C(T(2,3); j, 6j + 1)
\ \rightarrow \
k_{IX}(j).
\]
\end{lemma}
%%%%%%%%%%%%%%%%%%%%%%%%%%%%%%%%%%%%
%
%
%

The proof of Theorem~\ref{thm:main} is divided into two parts:
In the first half, starting with Baker's Dehn surgery description
in \cite{Ba, Ba3},
we study the knots by usual diagrams.
In the second half, we will use divide presentations.
We will introduce a convenient method,
which we call {\it Couture move},
to deform generalized L-shaped curves.
It was pointed out in the private communication 
of the author and Olivier Couture.
With Couture moves, the proof get much geometric,
intuitive and shorter.
The author's old proof of Theorem~\ref{thm:main} was troublesome
braid calculus.
We will show Lemma~\ref{lem:Pj}
in the case $j < -1$ by Couture moves, as a demonstration,
in Subsection~\ref{sbsec:couture}.

\medskip

On the relation between the surgery coefficient and 
the area of the region (of the curve), 
the formula in Theorem~1.4 in \cite{Y4} is generalized to:
%
%
%%%%%%%%%%%%%%%%%%%%%%%%%%%%%%%%%%%%
\begin{lemma}~\label{lem:area}
On the divide presentation of $k_{\mathcal{X}}(j)$ 
by the generalized L-shaped curve $P_{\mathcal{X}}(j)$
in Theorem~\ref{thm:main},
the area, the number of concave points (of the region) of $P_{\mathcal{X}}(j)$
and the coefficient of the lens space surgery
along $k_{\mathcal{X}}(j)$ satisfy
\[
[ \textrm{Area} - \sharp \{\textrm{Concave points} \}] 
- \textrm{Surgery coefficient}
 = 0 \textrm{ or } 1.
\]
\end{lemma}
%%%%%%%%%%%%%%%%%%%%%%%%%%%%%%%%%%%%
%
%
see Definition~\ref{def:concave}
for the precise definition of concave points.
Lemma~\ref{lem:area} will be verified by Table~\ref{tbl:area}
in Subsection~\ref{sbsec:spradC}.
%
%
%%%%%%%%%%%%%%%%%%%%%%%%%%%%%%%%%%%%
\begin{question}
Study divide knots $L(P)$ presented by generalized L-shaped curves
$P = X \cap \mathcal{L}$.
By Lemma~\ref{lem:area}, $c(\mathcal{L}) = 
\textrm{Area}(\mathcal{L}) - \sharp \{\textrm{Concave points of } \mathcal{L} \}$
can be an expected coefficient for exceptional Dehn surgery of $L(P)$.
Study $c(\mathcal{L})$ and $c(\mathcal{L})-1$ surgeries along $L(P)$.
\end{question}
%%%%%%%%%%%%%%%%%%%%%%%%%%%%%%%%%%%%
%
%

This paper is organized as follows. In the next section,
we review theory of A'Campo's divide knots and links briefly
and generalize L-shaped plane curve and decide 
the parametrizing notation.
In Section~\ref{sec:sporadic}, we review the construction of the knots $k_{\mathcal{X}}(j)$,
give a precise definition of the plane curves $P(j)$
and prove Theorem~\ref{thm:main} and the lemmas.
\par

%%%%%%%%%%%%%%%%%%%%%%%%%%%%%%%%%%%%%%%%%%
%%%                   Section 2
%%%%%%%%%%%%%%%%%%%%%%%%%%%%%%%%%%%%%%%%%%
\section{Divide knots and plane curves}~\label{sec:curve}
We review theory of A'Campo's divide knots and links briefly.
We are interested in plane curves constructed 
as intersection of the $\pi /4$-lattice $X$ and a region.
We define a generalized L-shaped plane curve and decide 
the parametrizing notation.
%%%%%%%%%%%%%%%%%%%%%%%%%%%%%%%%%%%%
%%%%%%%%%%%%%%%%%%%%%%%%%%%%%%%%%%%%
\subsection{Torus knots}~\label{sbsec:torus}
\begin{figure}[h]
\begin{center}
\includegraphics[scale=0.4]{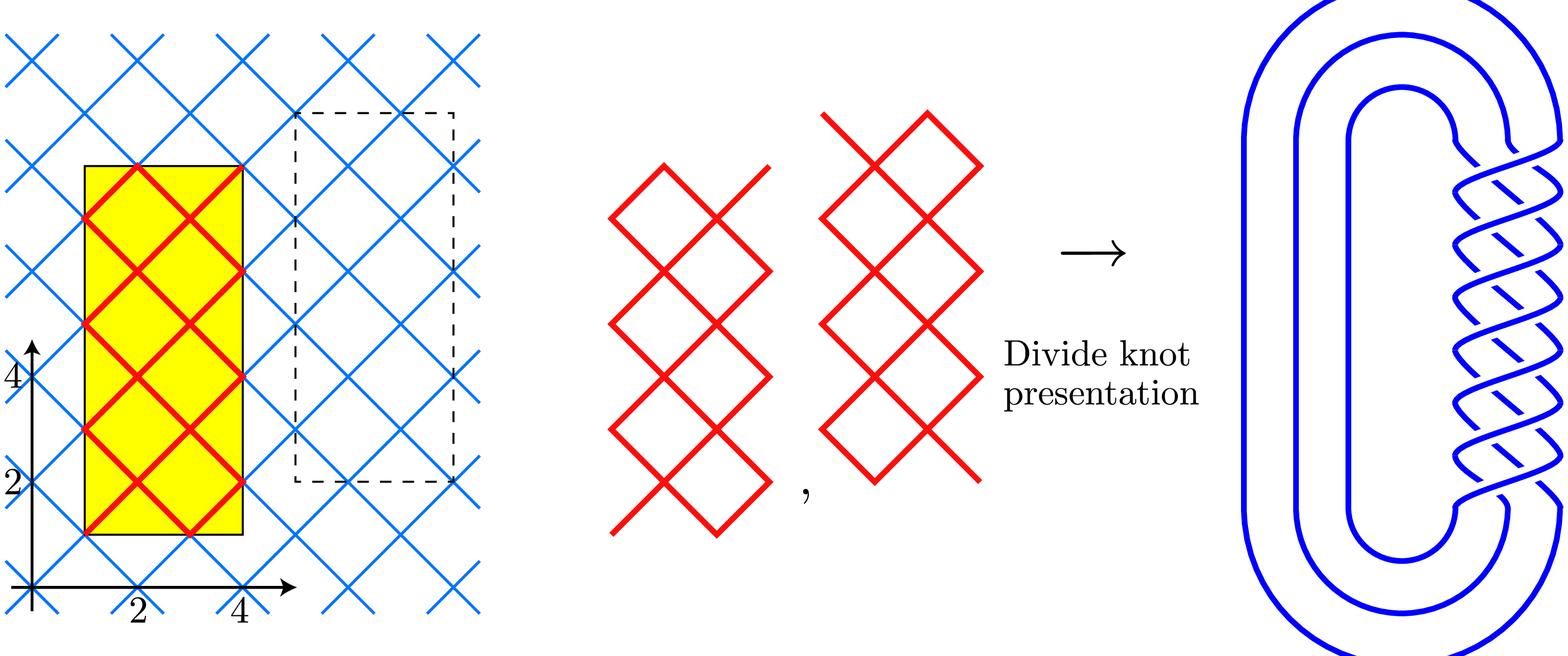}
\caption{A billiard curve $B(a,b)$ presents a torus knot $T(a,b)$ (ex. $T(3,7)$)}
\label{fig:Torus}
\end{center}
\end{figure}

We start with a presentation of a (positive) torus knot as a divide knot.
Let $(a,b)$ be a pair of positive integers
and $B(a,b)$ a curve defined as an intersection
of the $\pi /4$-lattice $X$
and an $a \times b$ rectangle $\mathcal{R}(a,b)$
whose every vertex is placed at a lattice point ($\in \mathbb{Z}^2$),
see Figure~\ref{fig:Torus}.
If $(a,b)$ is coprime, $B(a,b)$ is 
a billiard curve in $\mathcal{R}(a,b)$ with slope $\pm1$.
%
%
%%%%%%%%%%%%%%%%%
%%%%%%%%%%%%%%%%%
%
%
%%%%%%%%%%%%%%%%%%%%%%%%%%%%%%%%%%%%
\begin{lemma}~\label{lem:GHY}
{\rm (Couture--Perron \cite{CP}, Goda--Hirasawa--Y \cite{GHY}, see also \cite{AGV})}\\
The curve $B(a,b)$
presents the torus link $T(a,b)$ as a divide knot.
\end{lemma}
%%%%%%%%%%%%%%%%%%%%%%%%%%%%%%%%%%%%
%
%
Strictly, the curve $B(a,b)$ depends on the placement of the rectangle,
whether the left-bottom corner of the region 
is a terminal of the curve or not, see Figure~\ref{fig:Torus} again.
Even if $(a,b)$ is not coprime (i.e., case of a torus link),
the curve $B(a,b)$ in either choice presents $T(a,b)$ (\cite{GHY}).
If $(a,b)$ is coprime, the proof is easy:
the reflection along $x$-axis maps one to the other.
\par

%%%%%%%%%%%%%%%%%%%%%%%%%%%%%%%%%%%%
%%%%%%%%%%%%%%%%%%%%%%%%%%%%%%%%%%%%
\subsection{Basic facts on divide knots}~\label{sbsec:divide}
The theory of A'Campo's {\it divide knots and links}
comes from singularity theory of complex curves.
A {\it divide} $P$ is (originally)
a relative, generic immersion of a 1-manifold
in the unit disk $D$ in $\mathbb{R}^2$.
A'Campo \cite{A1, A2, A3, A4}
formulated the way to associate to each divide $P$ a link
$L(P)$ in $S^3$. We regard $S^3$ as 
\[
S^3 = \{
(u,v) \in D \times T_u D
\, | \,
\vert u \vert^2 + \vert v \vert^2 = 1
\}
\]
and the original construction is 
\[
L(P) = \{
(u,v) \in D \times T_u D
\, | \,
u \in P, v \in T_uP, 
\vert u \vert^2 + \vert v \vert^2 = 1
\} \subset S^3,
\]
where $T_uP$ is the subset consisting of vectors tangent to $P$ 
in the tangent space $T_u D$ of $D$ at $u$.
In the present paper, we regard a PL (piecewise linear) plane curve as a divide
by smoothing corners and controlling the size.

\begin{figure}[h]
\begin{center}
\includegraphics[scale=0.5]{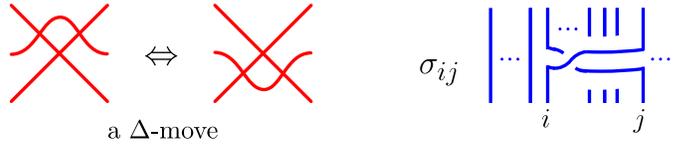} \par
\caption{Basics on divide knots}
\label{fig:div}
\end{center}
\end{figure}

Some characterizations of (general) divide knots and links are known,
and some topological invariants of $L(P)$
can be gotten from the divide $P$ directly.
Here, we list some of them.
%
%
%%%%%%%%%%%%%%%%%%%%%%%%%%%%%%%%%%%%%%
\begin{lemma}~\label{lem:divide}
{\rm  ((1)--(7) by A'Campo \cite{A2},
(8) by Hirasawa \cite{Hi}, Rudolph \cite{R})}
\begin{enumerate}
\item[(1)]
$L(P)$ is a knot (i.e., connected) if and only if
$P$ is an immersed arc.
\item[(2)]
If $L(P)$ is a knot, the unknotting number, the Seifert genus and 
the $4$-genus of $L(P)$ are all equal to
the number $d(P)$ of the double points of $P$.
\item[(3)]
If $P = P_1 \cup P_2$ is the image of an immersion of two arcs,
then the linking number of the two component link
$L(P) = L(P_1) \cup L(P_2)$ is equal to
the number of the intersection points between $P_1$ and $P_2$.
\item[(4)]
If $P$ is connected, then $L(P)$ is {\it fibered}.
\item[(5)]
Any divide link $L(P)$ is {\it strongly invertible}.
\item[(6)]
A divide $P$ and its mirror image $P!$ present
the same knot or link: $L(P!) =L(P)$.
\item[(7)]
If $P_1$ and $P_2$ are related by some $\Delta$-moves,
then the links $L(P_1)$ and $L(P_2)$ are isotopic:
If $P_1 \sim_{\Delta} P_2$ then $L(P_1) = L(P_2)$,
see Figure~\ref{fig:div}.
\item[(8)]
Any divide knot is a closure of a
{\it strongly quasi-positive} braid,  i.e.,
a product of some $\sigma_{ij}$ in Figure~\ref{fig:div}.
\end{enumerate}
\end{lemma}
%%%%%%%%%%%%%%%%%%%%%%%%%%%%%%%%%%%%
%
%
%%%%%%%%%%%%%%%%%
%%%%%%%%%%%%%%%%%
%
%
%For theory of divide knots, see also Rudolph's \lq\lq {\bf C}-link\rq\rq \ %
%in \cite{R} and \cite{Ch, HW}.
For theory of divide knots, see also \cite{Ch, HW}
and ^^ ^^ transverse $\mathbb{C}$-links" defined by Rudolph \cite{R}.
In \cite{CP} Couture and Perron pointed out a method
to get the braid presentation from the divide 
in a restricted cases, called \lq\lq ordered Morse\rq\rq \ %
divides.
We can apply their method.
It is a special case of Hirasawa's method in \cite{Hi}.
\begin{figure}[h]
\begin{center}
\includegraphics[scale=0.7]{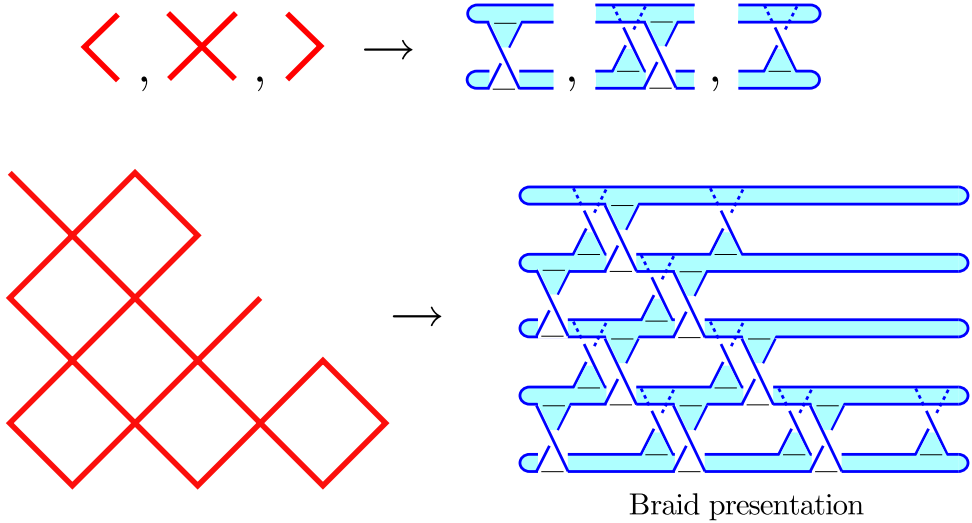}
\caption{Couture--Perron's method}
\label{fig:CP}
\end{center}
\end{figure}

Finally we recall an operation ^^ ^^ adding a square" 
on divides $P$ and its contribution to the divide links $L(P)$.
%
%
%%%%%%%%%%%%%%%%%%%%%%%%%%%%%%%%%%%%%%
\begin{lemma}[Lemma~4.2 in \cite{Y4}]~\label{lem:addS}%[Lemma~4.2]{Y4}
Adding a square on an L-shaped curve $P$ along an edge $l$ (of the region)
corresponds to a right handed full-twist on the divide knot 
$L(P)$ along the unknotted component defined by $l$.
\end{lemma}
%%%%%%%%%%%%%%%%%%%%%%%%%%%%%%%%%%%%%%
%
%
Adding a square is related to ^^ ^^ blow-down". 
Here, blow-down along $y$-axis ($x=0$) is
the coordinate transformation from $(x,y)$ to $(x', y')$
by $x'=x, y' = yx$ (ex. $y^2 = x+\epsilon$ become ${y'}^2 = {x'}^2 (x'+\epsilon)$)
%
%
%%%%%%%%%%%%%%%%%
\begin{figure}[h]
\begin{center}
\includegraphics[scale=0.3]{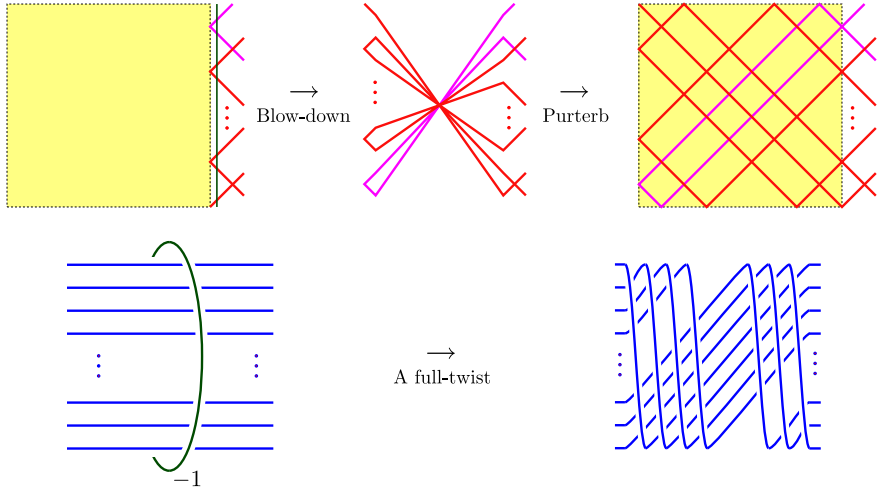}
\caption{Adding a square}
\label{fig:addS}
\end{center}
\end{figure}
%%%%%%%%%%%%%%%%%
%
%
\par 

%%%%%%%%%%%%%%%%%%%%%%%%%%%%%%%%%%%%
%%%%%%%%%%%%%%%%%%%%%%%%%%%%%%%%%%%%
\subsection{Curves defined by regions}~\label{sbsec:region}
In $xy$-plane, 
a lattice point or an integer vector $(k,l) \in \mathbb{Z}^2$ is called
{\it even} or {\it odd},
if $k+l$ is even or odd, respectively.
Double points of the $\pi /4$-lattice
$X = \{ (x,y) \vert \cos \pi x = \cos \pi y \}$ are at even points.
In this paper, we are interested in curves constructed 
as intersection of $X$ and a region, 
as a generalization of Lemma~\ref{lem:GHY}.

The lattice $X$ has some symmetries:
We let $r_x$ denote the reflection along $y$-axis, $r_x (x,y) =(-x,y)$,
$R_{\pi/2}$ the $\pi/2$-rotation (along the origin), and 
$+ \vec{v}_{\text{ev}}$ a parallel translation 
by an \underline{even} vector.
Symmetry of the lattice $X$ is generated by 
$r_x, R_{\pi/2}$ and some $\vec{v}_{\text{ev}}$.
A curve constructed 
as intersection of $X$ and a region $\mathcal{R}$
does not change by the action (on $\mathcal{R}$) of the symmetry of $X$.
We also use a parallel translation $+ \vec{v}_{\text{od}}$ by an odd vector
to place a region $\mathcal{R}$ well.

%
%
%%%%%%%%%%%%%%%%%%%%%%%%%%%%%%%%%%%%%%
\begin{definition}[Condition of regions]~\label{def:concave}
We are interested in curves constructed 
as intersection of $X$ and a region $\mathcal{R}$.
We formulate the conditions on regions:
\begin{enumerate}
\item[(i)] A region $\mathcal{R}$ is a union of a 
finite number of rectangles.
\item[(ii)] Each vertex of $\mathcal{R}$ is at a lattice point.
\item[(iii)] Each edge of the rectangles in $\mathcal{R}$ is 
horizontal or vertical.
\item[(iv)] 
Difference vectors of any pair of concave points 
of the region ${\mathcal R}$ are even.
\end{enumerate}
Here, a concave point in (iv) is defined as follows:
A boundary point $p$ of a region $\mathcal{R}$
is called a {\it concave point} (of the region) if 
a neighborhood of $\mathcal{R}$ at $p$ is locally 
homeomorphic to that of $\{ x \leq 0\} \cup \{ y \leq 0\}$ at $(0,0)$ by
the symmetry of $xy$-plane, generated by 
$r_x$, $R_{\pi/2}$ and $+ \vec{v}$ (by an even or an odd integer vector),
see Figure~\ref{fig:gLreg}.
\end{definition}
%%%%%%%%%%%%%%%%%%%%%%%%%%%%%%%%%%%%%%
%
%
If a concave point $p$ of $\mathcal{R}$ is at even point, 
then the curve $X \cap \mathcal{R}$ is not generic at $p$,
i.e., a terminal point overlaps with an interior point of the curve.
We are concerned only with a generic immersed curves.
By the condition (iv), all concave points of either
$X \cap \mathcal{R}$ or $X \cap (\mathcal{R} + \vec{v}_{\text{od}})$
are placed at even points, and it defines a generic immersed curve.
%
%
%%%%%%%%%%%%%%%%%%%%%%%%%%%%%%%%%%%%%%
\begin{definition}~\label{def:RC}
For a region $\mathcal{R}$ satisfying the condition (i),(ii),(iii) and (iv),
either $X \cap \mathcal{R}$ or 
$X \cap (\mathcal{R} + \vec{v}_{\text{od}})$ is a generic immersed curve.
We choose the generic one and define it 
as {\it a curve defined by the region $\mathcal{R}$},
see an example
($X \cap (\mathcal{R} + \vec{v}_{\text{od}})$ is chosen) in Figure~\ref{fig:gLreg}.
\end{definition}
%%%%%%%%%%%%%%%%%%%%%%%%%%%%%%%%%%%%%%
%
%
We describe a curve 
$X \cap \mathcal{R}$ by describing (and parametrizing) 
the region $\mathcal{R}$ by using $xy$-coordinates.
\par

%%%%%%%%%%%%%%%%%%%%%%%%%%%%%%%%%%%%
%%%%%%%%%%%%%%%%%%%%%%%%%%%%%%%%%%%%
\subsection{Generalized L-shaped curves}~\label{sbsec:gL}
We define {\it generalized L-shaped curves}.
It is an extension of ^^ ^^ L-shaped curves" in \cite[Section~3.2]{Y4},
but the notation (parametrization) is changed.
%
%
%%%%%%%%%%%%%%%%%%%%%%%%%%%%%%%%%%%%%%
\begin{definition}[Generalized L-shaped region at the origin]~\label{def:gL}
See Figure~\ref{fig:gLreg}.\\
Let $n$ be a positive integer with $n > 1$. 
We let 
\begin{eqnarray}~\label{eq:param}
[ (a_i,b_i) ] := 
[
(a_1, b_1), (a_2,b_2), \cdots, (a_n,b_n)
]
\end{eqnarray}
denote a sequence of lattice points ($\in \mathbb{Z}^2$)
in $xy$-plane satisfying 
\begin{center}
$0 < a_1 < a_2 < \cdots < a_n$ and 
$b_1 > b_2 > \cdots > b_n > 0$.
\end{center}
We define a region $R[ (a_i,b_i) ]$ in $xy$-plane by 
\begin{eqnarray}~\label{eq:gL}
R[ (a_i,b_i) ]
= 
\bigcup_{i=1}^n
\{
(x,y)
\vert \,
0 \leq x \leq a_i,
0 \leq y \leq b_i
\}.
\end{eqnarray}
%
%
%%%%%%%%%%%%%%%%%
\begin{figure}[h]
\begin{center}
\includegraphics[scale=0.4]{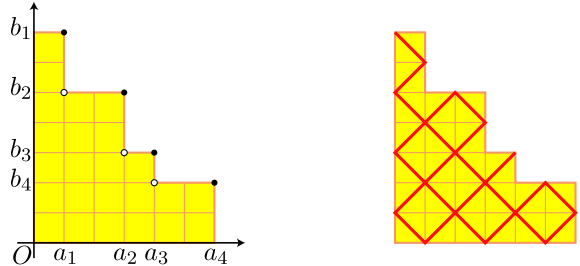}
\caption{Type $R[(1,7),(3,5),(4,3),(6,2)]$ ($=P(3)$)}
\label{fig:gLreg}
\end{center}
\end{figure}
%%%%%%%%%%%%%%%%%
%
%
We will call this region
{\it a generalized L-shaped region of type $[ (a_i,b_i) ]$
of length $n$} (at the origin).

If such a region defines a generic immersed curve
in the sense of Definition~\ref{def:RC},
we call the curve
{\it generalized L-shaped curve of type $[ (a_i,b_i) ]$}.
\end{definition}
%%%%%%%%%%%%%%%%%%%%%%%%%%%%%%%%%%%%%%
%
%
A generalized L-shaped region of type $[ (a_i,b_i) ]$ of length $n$
has $n-1$ concave points 
at the coordinate $(a_i,b_{i+1})$ with $1 \leq i \leq n-1$.
Note that the L-shaped region $L[a_1, a_2; b_1, b_2]$ defined in \cite[Definition~3.3]{Y4}
is redefined $R[(a_1,b_2),(a_2,b_1)]$ of length $2$ here.

It is easy to see:
%
%
%%%%%%%%%%%%%%%%%%%%%%%%%%%%%%%%%%%%%%
\begin{lemma}~\label{lem:gLarea}
The area of a generalized L-shaped region of type $[ (a_i,b_i) ]$
of length $n$ is
\[
\text{Area }(R[ (a_i,b_i) ])
= 
\sum_{i=1}^n a_i b_i - \sum_{i=1}^{n-1} a_{i} b_{i+1}.
\]
\end{lemma}
%%%%%%%%%%%%%%%%%%%%%%%%%%%%%%%%%%%%
%
%
%%%%%%%%%%%%%%%%%%%%%%%%%%%%%%%%%%%%
\begin{question}
Find a formula on the numbers of circle and arc components of 
(generic) generalized L-shaped curves.
When is a generalized L-shaped curve a generic immersed arc?
\end{question}
%%%%%%%%%%%%%%%%%%%%%%%%%%%%%%%%%%%%
%
%

%%%%%%%%%%%%%%%%%%%%%%%%%%%%%%%%%%%%
%%%%%%%%%%%%%%%%%%%%%%%%%%%%%%%%%%%%
\subsection{Couture move}~\label{sbsec:couture}
\begin{figure}[h]
\begin{center}
\includegraphics[scale=0.4]{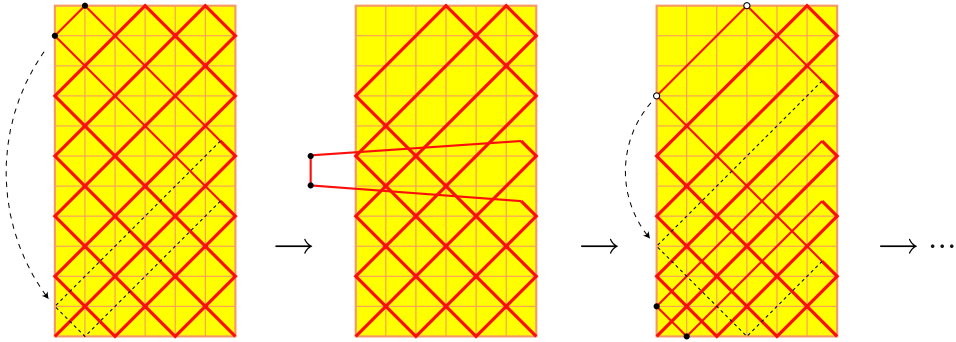}
\caption{Couture move consists of some $\Delta$-moves}
\label{fig:PfCmv}
\end{center}
\end{figure}

We introduce a convenient method {\it Couture move}.
It was pointed out in the private communication 
of the author and Olivier Couture in the opportunity of 
a conference
^^ ^^ Singularities, knots, and mapping class groups in memory of Bernard Perron"
held in Sept.\ 2010.
The purpose was to present torus knots by generalized L-shaped curves
(other than $B(a,b)$ in Lemma~\ref{lem:GHY}) as divide knots.
Here we characterize the move as follows:
%
%
%%%%%%%%%%%%%%%%%%%%%%%%%%%%%%%%%%%%%%
\begin{definition}~\label{def:Cmv}
We say that a deformation of plane curves (divides) is a {\it Couture move}
if 
\begin{enumerate}
\item[(1)]
It is from a curve $X \cap \mathcal{R}$
defined by a region $\mathcal{R}$ in the sense of Definition~\ref{def:RC},
\item[(2)]
The deformation consists of some $\Delta$-moves, and
\item[(3)]
The resulting curve is a curve $X \cap \mathcal{R}'$
defined by another region $\mathcal{R}'$.
\end{enumerate}
\end{definition}
%%%%%%%%%%%%%%%%%%%%%%%%%%%%%%%%%%%%%%
%
%
As a demonstration, we prove Lemma~\ref{lem:Pj} in the case $j < -1$.
%
%
%%%%%%%%%%%%%%%%%
\begin{figure}[h]
\begin{center}
\includegraphics[scale=0.4]{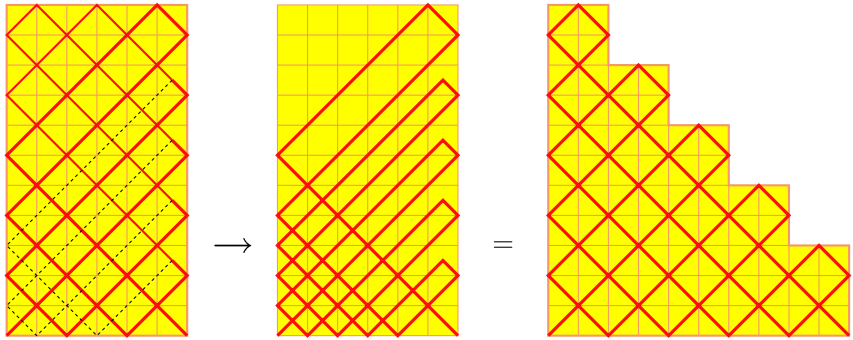}
\caption{Couture move (ex. $T(6,11)$)}
\label{fig:Cmv2}
\end{center}
\end{figure}
%%%%%%%%%%%%%%%%%
%
%
%%%%%%%%%%%%%%%%%%%%%%%%%%%%%%%%%%%%%%
\begin{lemma}~\label{lem:Tn2n-1}
Let $n$ be an integer with $n>1$.
We define a sequence $[(a_i,b_i)^{(n)}]$ of lattice points 
of length $n-1$ by
\[
(a_i, b_i)^{(n)} = (2i, 2n+1-2i) \quad
(i=1, 2, ..., n-1).
\]
The generalized L-shaped curve
of type $[(a_i,b_i)^{(n)}]$
presents the torus knot $T(n,2n-1)$ as a divide knot.
\end{lemma}
%%%%%%%%%%%%%%%%%%%%%%%%%%%%%%%%%%%%
{\bf Proof}. The method is shown in Figure~\ref{fig:Cmv2}, which is 
an example from $B(5,11)$ to the generalized L-shaped curve
of type $[(2,11), (4,9), (6,7), (8,5), (10,3)]$.
It consists of some $\Delta$-moves, see Figure~\ref{fig:PfCmv}.
%
%
%%%%%%%%%%%%%%%%%
%%%%%%%%%%%%%%%%%
%
%
\qed
\par

%%%%%%%%%%%%%%%%%%%%%%%%%%%%%%%%%%%%%%%%%%
%%%                   Section 3
%%%%%%%%%%%%%%%%%%%%%%%%%%%%%%%%%%%%%%%%%%
\section{Details on sporadic knots and Proof}~\label{sec:sporadic}
We give a precise definition of the plane curves $P(j)$,
verify the formula in Lemma~\ref{lem:area} and prove Theorem~\ref{thm:main} and the lemmas.
We start the proof with Baker's description in \cite{Ba, Ba3} of sporadic knots,
and we use Couture moves on divides in the second half of the proof.

%%%%%%%%%%%%%%%%%%%%%%%%%%%%%%%%%%%%
%%%%%%%%%%%%%%%%%%%%%%%%%%%%%%%%%%%%
\subsection{Precise definition of Curves}~\label{sbsec:spradC}
We define divides $P(j)$ by the method introduced in the last section.
%
%
%%%%%%%%%%%%%%%%%%%%%%%%%%%%%%%%%%%%%%
\begin{definition}[Precise definition of $P(j)$]~\label{def:p(j)}
For an integer $j$ ($j \not=0,-1$),
we define a sequence $[(a_i,b_i)^{(j)}]$ of lattice points
$(a_i,b_i)^{(j)}$ as follows. 
\par
\medskip

(Case $j > 0$) \ 
Starting with $[(a_i,b_i)^{(1)}] = [(1,3),(2,1)]$, 
we define $[(a_i,b_i)^{(j)}]$ inductively with respective to $j$, 
by
\[
\begin{cases}
(a_1, b_1)^{(j)} = (1, 2j+1) & (i=1) \\
(a_i, b_i)^{(j)} = (b_{j+2-i}, a_{j+2-i})^{(j-1)} + (1,1) 
& (1< i \leq j+1) \\
\end{cases}
\]

(Case $j < -1$) \ 
We define $[(a_i,b_i)^{(j)}]$ by
\[
(a_i, b_i)^{(j)} = (2i, -2j +1 -2i)  
\quad (1\leq i \leq - j)
\]
We define a plane curve $P(j)$ as a generalized L-shaped curve
of type $[(a_i,b_i)^{(j)}]$,
whose length is $j+1$ (if $j>0$) or $-j$ (if $j<-1$).
See and verify the examples in Figure~\ref{fig:SporadPN}.
\end{definition}
%%%%%%%%%%%%%%%%%%%%%%%%%%%%%%%%%%%%%%
%
%
By Lemma~\ref{lem:gLarea}, it is easy to see 
%
%
%%%%%%%%%%%%%%%%%%%%%%%%%%%%%%%%%%%%%%
\begin{lemma}~\label{lem:areaPj}
\[
Area ( R[(a_i,b_i)^{(j)}])
= 
\begin{cases}
2j^2 + 2j  & (j > 0) \\
2j^2 & (j < -1) \\
\end{cases}.
\]
\end{lemma}
%%%%%%%%%%%%%%%%%%%%%%%%%%%%%%%%%%%%
%
%
The plane curve $P_{\mathcal{X}}(j)$ is constructed 
by adding a square twice as in Definition~\ref{def:PXj}.
Since a square addition along an edge $E$ of length $l(E)$
increases the area by $l(E)^2$, 
the area of the region of $P_{IX}(j)$ with $j>0$ is
calculated as 
\begin{eqnarray*}
Area ( R[(a_i,b_i)^{(j)}]) 
+ l(E_a)^2
+ l(\overline{E_b})^2
& = &
(2j^2 + 2j) + (2j)^2 + (4j+1)^2 \\
& = & 22 j^2 +10j + 1
\end{eqnarray*}
On the other hand, 
the area of the region of $P_{X}(j)$ with $j>0$ is
calculated as
\begin{eqnarray*}
Area ( R[(a_i,b_i)^{(j)}]) 
+ l(E_b)^2
+ l(\overline{E_a})^2
& = &
(2j^2 + 2j) + (2j+1)^2 + (4j+1)^2 \\
& = & 22 j^2 +14j + 2
\end{eqnarray*}
We calculate them also in the cases $j <-1$ and list them in Table~\ref{tbl:area}.
Lemma~\ref{lem:area} is proved by Table~\ref{tbl:area}.
%
%
%%%%%%%%%%%%%%%%%
\begin{table}[h]
\begin{center}
\begin{tabular}{|l|c|c|c|c|}
\hline
Curve  & Coeff. $p$ of $k_{\mathcal{X}}(j)$
& Area (Case $j>0$) & Area (Case $j< -1$) \\
\hline
$P_{IX}(j)$ & $22 j^2 + 9j +1$ 
& $22 j^2 + 10j +1$ & $22 j^2 +8j +1$ \\
\hline
$P_{X}(j)$ & $22 j^2 + 13j +2$ & 
$22 j^2 + 14j +2$ & $22 j^2 +12j +2$ 
\\
\hline
\end{tabular}
\caption{Area of the plane curve $P_{\mathcal{X}}(j)$}
\label{tbl:area}
\end{center}
\end{table}
%%%%%%%%%%%%%%%%%
%
%

Next, we calculate and verify that the numbers of double points 
of $P(j)$ and $P_{\mathcal{X}}(j)$.
By Lemma~\ref{lem:divide}(2),
they are equal to the genus of the presented knots
$L(P(j))$ and $L(P_{\mathcal{X}}(j)) = k_{\mathcal{X}}(j)$, respectively.
Since $L(P(j))$ is $T(j, 2j+1)$ by Lemma~\ref{lem:Pj}, we have:
%
%
%%%%%%%%%%%%%%%%%%%%%%%%%%%%%%%%%%%%%%
\begin{lemma}~\label{lem:dPj}
The number $d(P(j))$ of double points of the plane curve $P(j)$ is
equal to the genus of the torus knot $T(j, 2j+1)$:
\[
d( P(j) )
= 
\frac{(\vert j \vert -1)(\vert 2j+1 \vert -1)}2
= 
\begin{cases}
j (j-1) & (j > 0) \\
(j+1)^2 & (j < -1) \\
\end{cases}.
\]
\end{lemma}
%%%%%%%%%%%%%%%%%%%%%%%%%%%%%%%%%%%%
%
%
Since the number of double points increases by $l(l-1)/2$
by adding a square along an edge of length $l$,
$d( P_{IX}(j) )$ is calculated as 

(Case $j>0$)
\begin{eqnarray*}
d( P_{IX}(j) )
& = & 
d( P(j) ) 
+ \frac{l(E_a) \cdot (l(E_a) -1)}2
+ \frac{l(\overline{E_b}) \cdot ( l(\overline{E_b}) -1 )}2
\\
& = & 
j(j -1) + \frac{2j \cdot (2j-1)}2 +  \frac{(4j+1) \cdot 4j}2
\\
& = & 11j^2 
\end{eqnarray*}

(Case $j<-1$)
\begin{eqnarray*}
d( P_{IX}(j) )
& = & 
(j +1)^2 + \frac{(-2j) \cdot (-2j-1)}2 +  \frac{(-4j-1) \cdot (-4j-2)}2
\\
& = & 11j^2 + 9j +2
\end{eqnarray*}
They are equal to the genus of the knots $k_{IX}(j)$, 
see Table~\ref{tbl:sporad}.
 We leave the case of TypeX to the readers.

%%%%%%%%%%%%%%%%%%%%%%%%%%%%%%%%%%%%
%%%%%%%%%%%%%%%%%%%%%%%%%%%%%%%%%%%%
\subsection{Proof of the main theorem}~\label{sbsec:torus}
In the first half of the proof, we study the knots in the usual diagram
and Dehn surgery description. 
In the second half, we will use divide presentation.

\medskip

\begin{figure}[h]
\begin{center}
\includegraphics[scale=0.4]{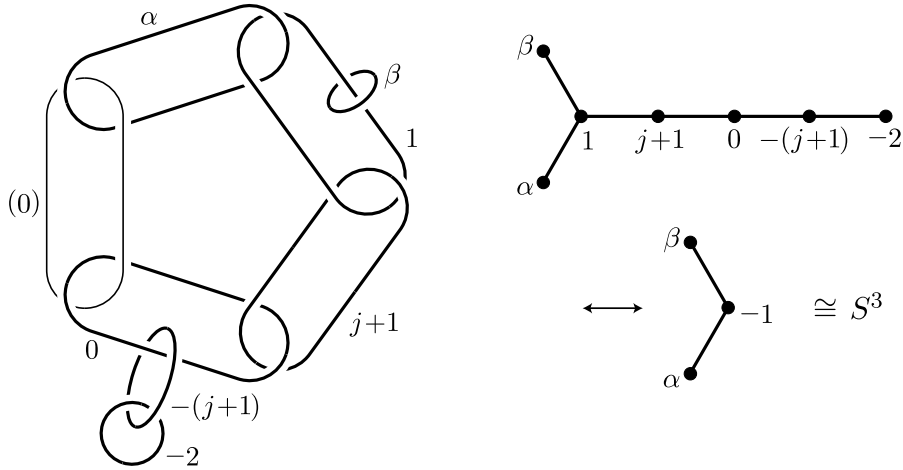}
\caption{Baker's description}
\label{fig:Baker}
\end{center}
\end{figure}

We start with Baker's Dehn surgery description of the knots 
$k_{\mathcal{X}}(j)$ in Figure~\ref{fig:Baker} from \cite{Ba}.
Throughout the paper, we fix 
\[
(\alpha, \beta ) = 
\begin{cases}
(-2,-3) & \text{ for TypeIX} \\
(-3,-2) & \text{ for TypeX}
\end{cases}.
\]
%
%
%%%%%%%%%%%%%%%%%
%%%%%%%%%%%%%%%%%
%
%
It is easy to see that 
the framed sublink of thick seven components presents $S^3$, 
by usual weighted tree diagram in the right half of Figure~\ref{fig:Baker}.
We call the sublink the non-trivial diagram of $S^3$.
In the resulting $S^3$, the thin component 
is the sporadic knot $k_{\mathcal{X}}(j)$ as a knot.

\medskip

\begin{figure}[h]
\begin{center}
\includegraphics[scale=0.4]{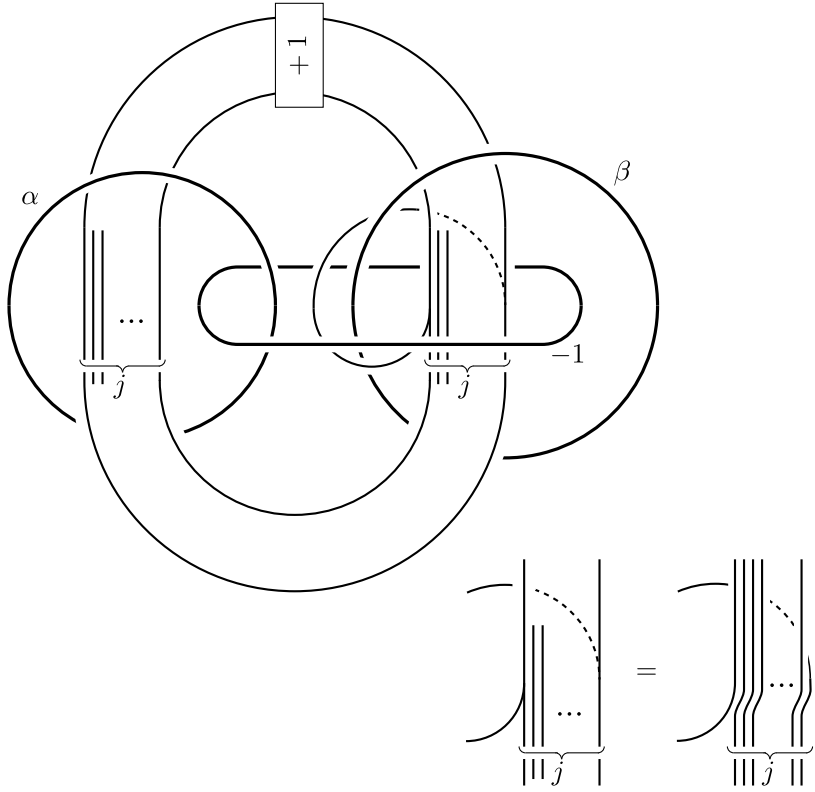}
\caption{Construction of $k_{\mathcal{X}}(j)$ : Link $L(j)$}
\label{fig:fromT}
\end{center}
\end{figure}

[Case $j>0$] \ First, we assume that $j>0$ until the final paragraph.
To get a usual diagram of the knot $k_{\mathcal{X}}(j)$,
we have to chase the thin component (with framing)
during the deformation from the non-trivial diagram to the empty diagram
of $S^3$. 
It is not easy but straight forward. In the middle of the process, 
we reach a diagram in Figure~\ref{fig:fromT}, 
where boxed $+1$ in the diagram means a right handed full-twist. 
The thin component $k(j)$ is a torus knot $T(j,j+1)$ with framing $j(j+1)$,
which is a coefficient of reducible Dehn surgery.
We name the four component link 
$L(j) = k(j)\cup u_{-1} \cup u_{\alpha} \cup u_{\beta}$,
where $u_x$ the $x$-framed unknotted component for $x= -1, \alpha, \beta$.

%
%
%%%%%%%%%%%%%%%%%
%%%%%%%%%%%%%%%%%
%
%

\begin{figure}[h]
\begin{center}
\includegraphics[scale=0.4]{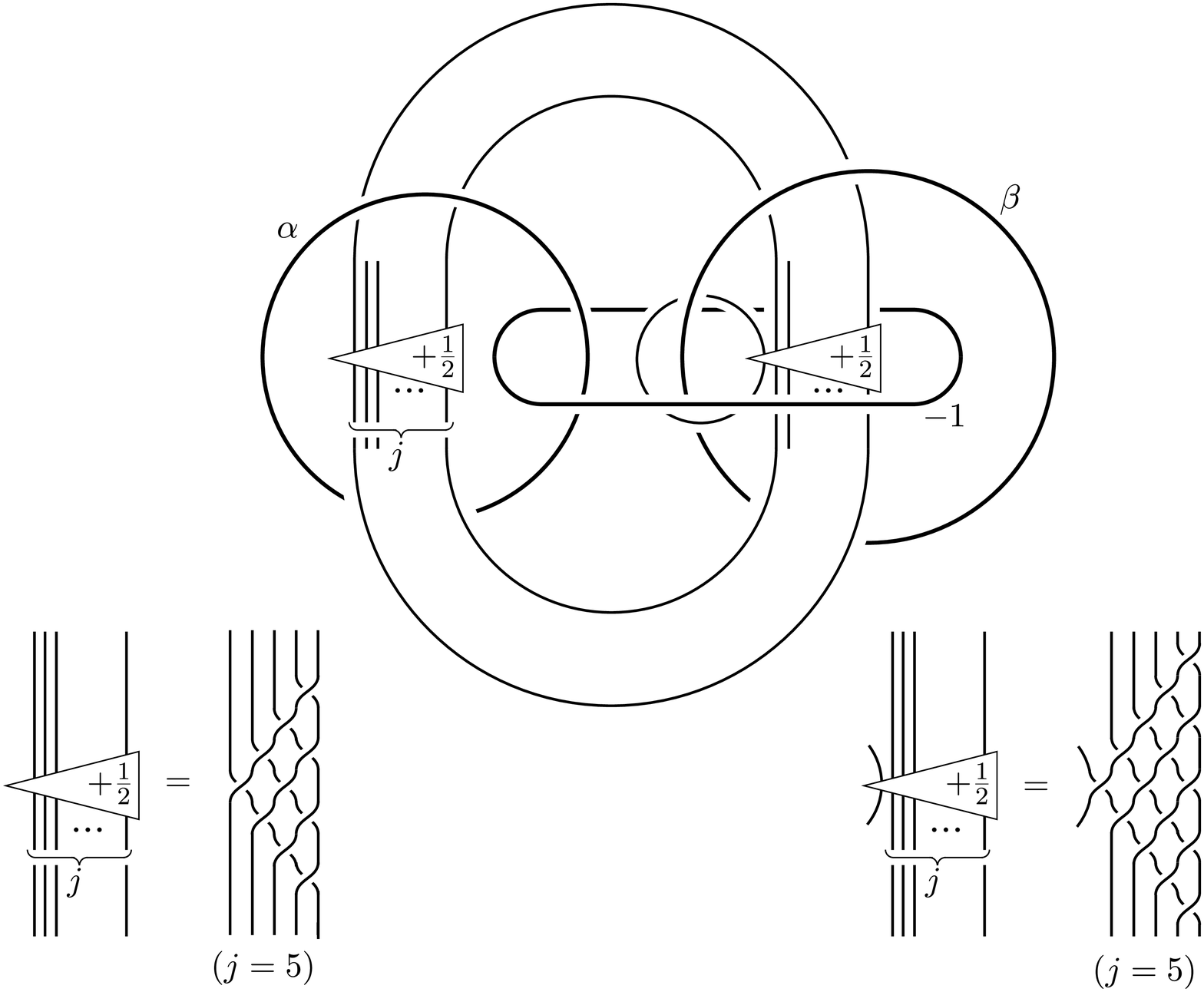}
\caption{Link $L(j)$ with a strong involution}
\label{fig:fromT2}
\end{center}
\end{figure}

\begin{figure}[h]
\begin{center}
\includegraphics[scale=0.4]{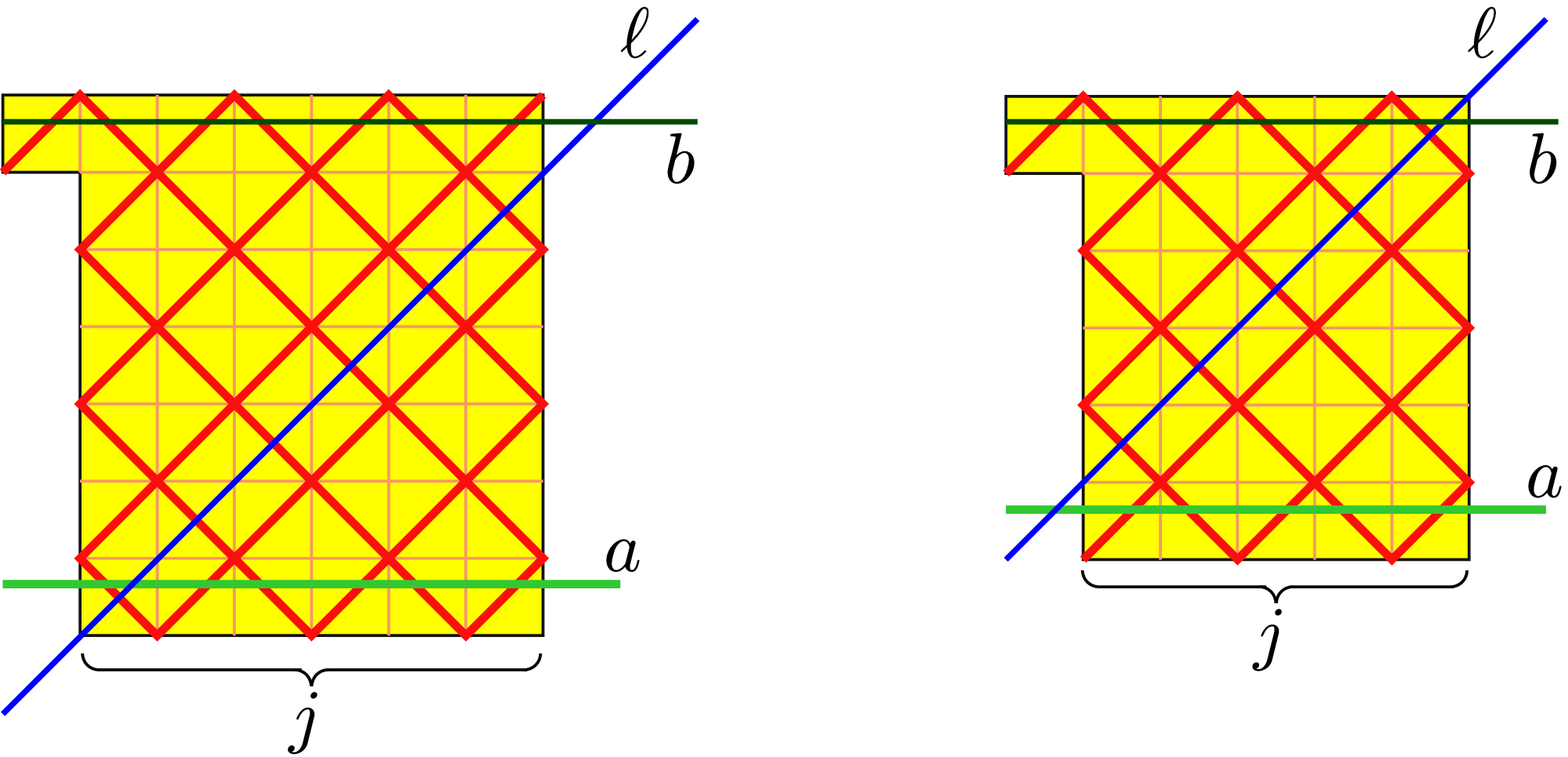} 
\caption{Divide presentation $C(j)$ of $L(j)$
\ (ex. $j =6$ and $j=5$)
}
\label{fig:DivEO}
\end{center}
\end{figure}

Second, we decompose the full-twist at the boxed $+1$
as two half-twists, denoted by $+ \frac12$ in triangles,
and deform the diagram by isotopy as in Figure~\ref{fig:fromT2}.
Note that the half-twist in the right is of $j+1$ strings. 
We can see that $L(j)$ is strongly invertible with respect to the horizontal axis,
see Lemma~\ref{lem:divide}(5).
%
%
%%%%%%%%%%%%%%%%%
%%%%%%%%%%%%%%%%%
%
%
As a quotient of the involution, 
ignoring the crossing data (over or under),
we have a plane curve properly immersed in the half plane.
The curve can be modified as in Figure~\ref{fig:DivEO}.
This is a divide presentation of $L(j)$. 
It can be checked by Couture--Perron's method in Figure~\ref{fig:CP}.
This process is related to the original construction of divide knots:
For a link of singularity of a complex plane curve, 
the divide is a real part of a ^^ ^^ good" perturbation 
(called real Morsification) of the equation of the singularity.

In the divide presentation of $L(j)$ in Figure~\ref{fig:DivEO},
a line segment $\ell$ is placed slightly 
different whether $j$ is even or odd.
We name the plane curve $C(j) = c(j) \cup \ell \cup a \cup b$,
where $c(j)$ presents $k(j) = T(j, j+1)$ by Lemma~\ref{lem:GHY},
since $c(j)$ is a generalized L-shaped curve of 
type $[(j,j+1),(j+1,1)]$ and is isotopic to $B(j,j+1)$.
The line segments $\ell, a, b$ presents 
$u_{-1}, u_{\alpha}, u_{\beta}$, respectively, 
as a divide presentation.
%
%
%%%%%%%%%%%%%%%%%
%%%%%%%%%%%%%%%%%
%
%
The linking matrix of $L(j)$ with a suitable orientation of the link
\[
\begin{pmatrix}
j(j+1) & j & j & j+1 \\
j & -1 & 1 & 1 \\
j &  1 & \alpha & 1 \\
j+1 &  1 & 1 & \beta \\
\end{pmatrix},
\]
is equal to the matrix of the number of intersection points
of the components of the divide by Lemma~\ref{lem:divide}(3).

\medskip

From now on, we go into the second half of the proof,
and study the knots and links by divide presentation.
We use $\Delta$-moves on divides freely, see Lemma~\ref{lem:divide}(7).
We have two (or three) steps:
(i) Blow-down along $\ell$ (take a full-twist along $u_{-1}$,
(ii) (Only if $j$ is odd) Modify the curve by some $\Delta$-moves,
and 
(iii) Deform the curve by Couture moves.
Our goal is the divide $P(j)$ with edges $\overline{a}, \overline{b}$,
where $\overline{a}$ (and $\overline{b}$) is a small parallel push-off of 
the bottom edge $E_a$ (the left edge $E_b$) of the region of $P(j)$
into the interior.

\begin{figure}[h]
\begin{center}
\includegraphics[scale=0.4]{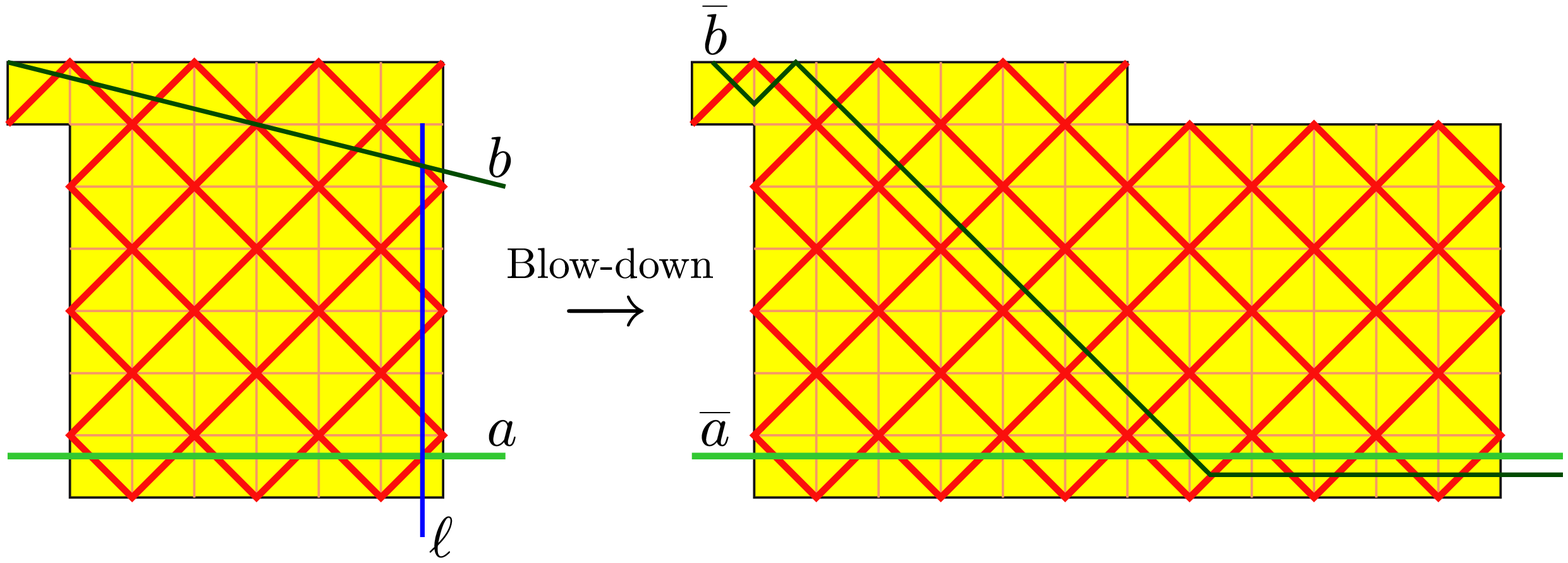} \\
\bigskip
\includegraphics[scale=0.4]{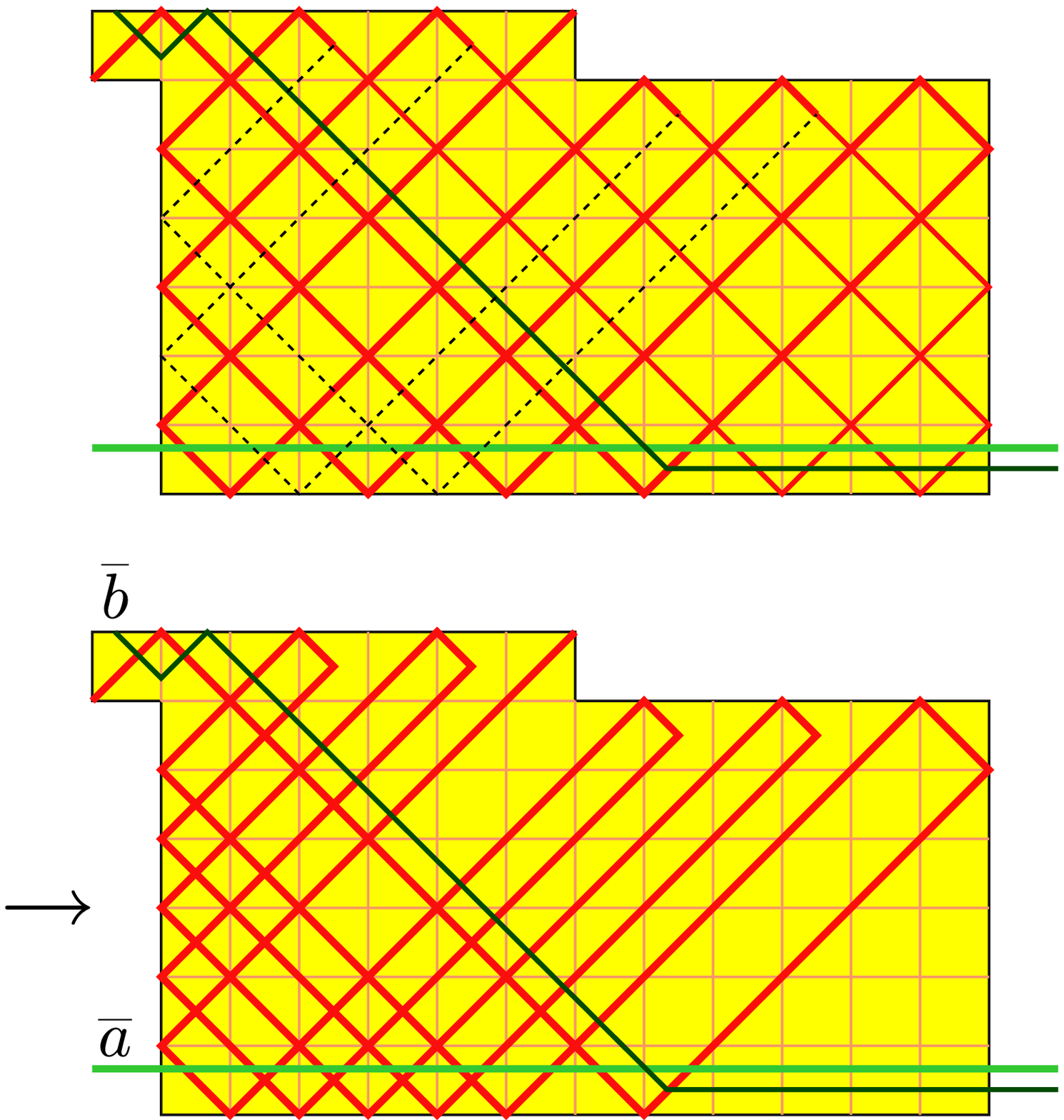}
\caption{Deformation: Case $j$ is even \ (ex. $j=6$)}
\label{fig:DivEv}
\end{center}
\end{figure}

\begin{figure}[h]
\begin{center}
\includegraphics[scale=0.4]{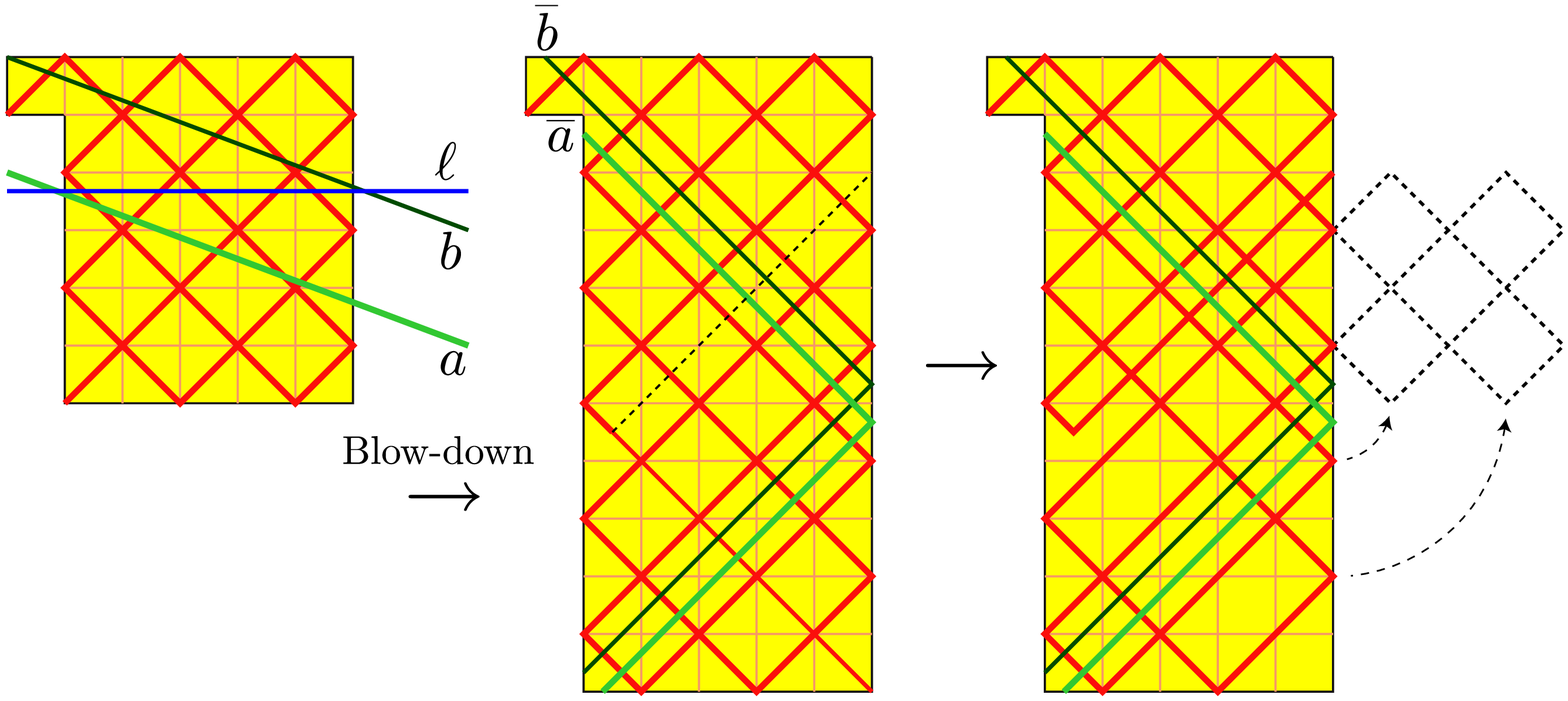} \\
\bigskip
\includegraphics[scale=0.4]{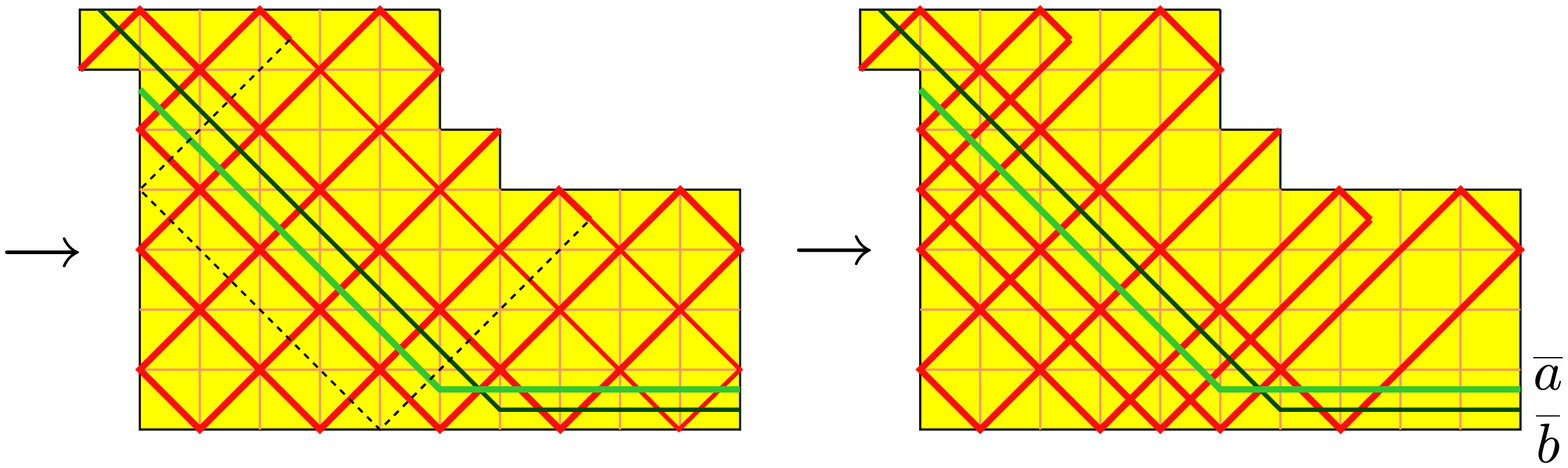}
\caption{Deformation: Case $j$ is odd \ (ex. $j=5$)}
\label{fig:DivOd}
\end{center}
\end{figure}

(Step (i)) \ We take a right handed full-twist of 
$k(j) \cup u_{\alpha} \cup u_{\beta}$ along the unknot $u_{-1}$.
We name the resulting link 
\[
\overline{L(j)} = 
\overline{k(j)} 
\cup \overline{u_{\alpha}}
\cup \overline{u_{\beta}}.
\]
This full-twist is done as a blow-down along the line $\ell$, 
i.e., by adding a square by Lemma~\ref{lem:addS}. 
In the case where $j$ is even, we slide $\ell$ and $b$ by $\Delta$-moves
as the first picture in Figure~\ref{fig:DivEv} and add a square. 
Otherwise, we add a square along the bottom edge as in Figure~\ref{fig:DivOd}.
%
%
%%%%%%%%%%%%%%%%%
%%%%%%%%%%%%%%%%%
%
%
%
%
%%%%%%%%%%%%%%%%%
%%%%%%%%%%%%%%%%%
%
%

(Step (ii))\ 
If $j$ is odd, the curve has a terminal point at the right bottom corner.
We move the terminal point (and its segment) up 
along the right edge,
as the second deformation of in Figure~\ref{fig:DivOd}.
We also slide the other added part at the bottom to the right 
by some $\Delta$-moves.

(Step (iii))\ 
The resulting curve is near L-shaped, but the lines
are not in the required position.
%$\overline{a}$ is not horizontal, $\overline{b}$ is not vertical yet.
Here we use Couture move,
see the second halves of 
Figure~\ref{fig:DivEv} in the case $j$ is even,
and 
Figure~\ref{fig:DivOd} in the case $j$ is odd.
By some obvious $\Delta$-moves,
we have the required curve 
$P(j) \cup \overline{a} \cup \overline{b}$,
which presents $\overline{L(j)}$ as a divide link.

The sublink $\overline{u_{\alpha}} \cup \overline{u_{\beta}}$
is a Hopf link and the framings are $(\alpha +1, \beta +1)$
$= (-1,-2)$ for TypeIX, or $= (-2,-1)$ for TypeX.
By the construction of $P_{\mathcal{X}}(j)$
in Definition~\ref{def:PXj}, and by the correspondence between
adding a square to a divide and a full-twist of the divide link 
in Lemma~\ref{lem:addS}, 
we have the required divide presentation 
of $k_{\mathcal{X}}(j)$ with $j>0$.

\bigskip

\begin{figure}[h]
\begin{center}
\includegraphics[scale=0.4]{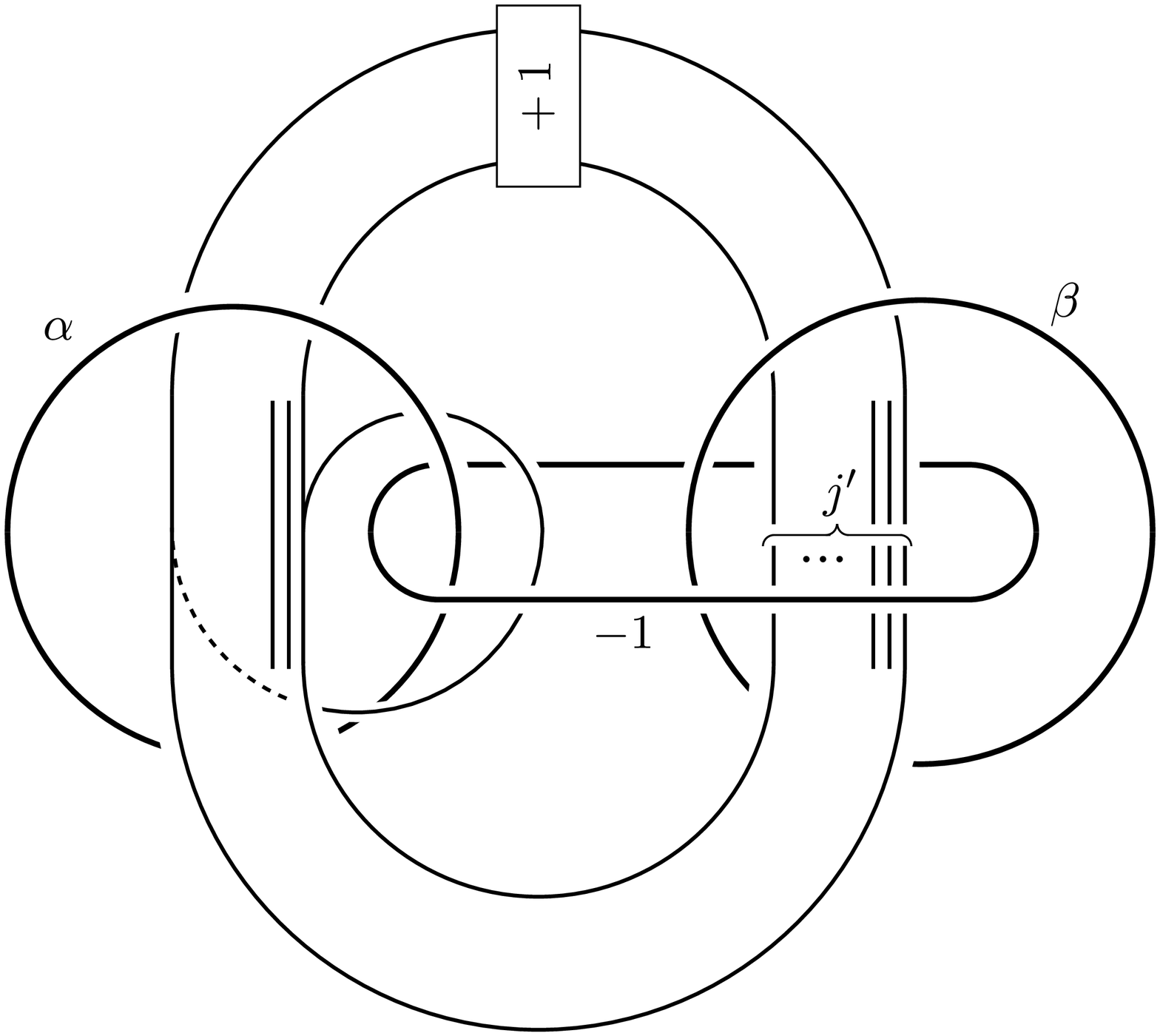}
\caption{Link $L(j)$ in the case $j<-1$ \ ($j=-(j'+1)$)}
\label{fig:fromTn}
\end{center}
\end{figure}

\begin{figure}[h]
\begin{center}
\includegraphics[scale=0.4]{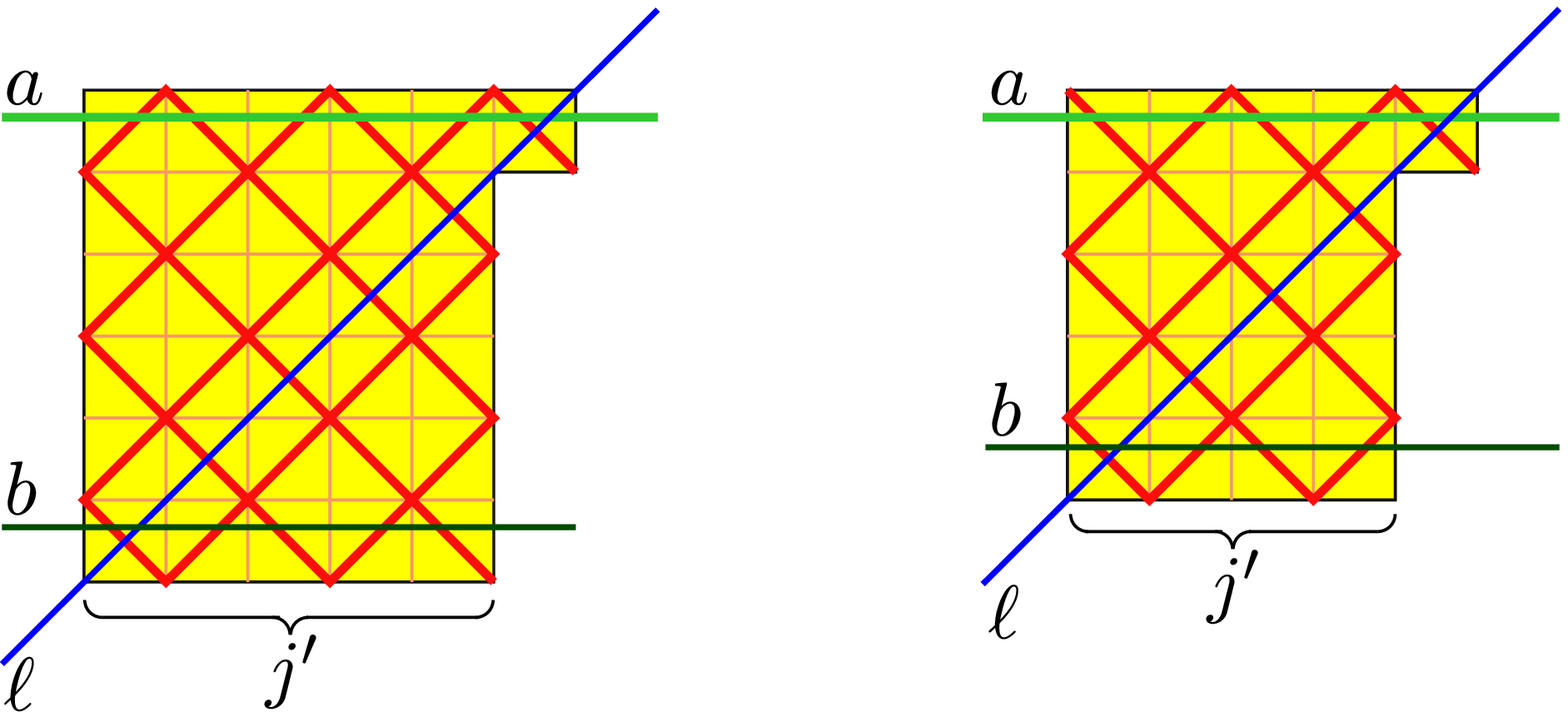} 
\caption{Divide presentation of $L(j)$ with $j<-1$
\ (ex. $j =-6$ and $j=-5$)
}
\label{fig:DivEOneg}
\end{center}
\end{figure}

\begin{figure}[h]
\begin{center}
\includegraphics[scale=0.4]{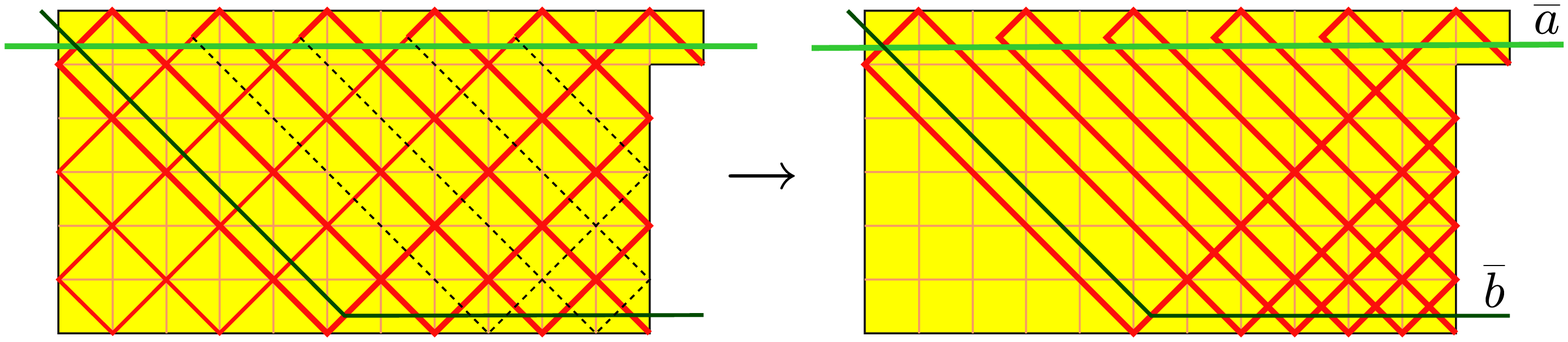}
\caption{Deformation: Case $j$ is even \ (ex. $j=-6$)}
\label{fig:DivEvneg}
\end{center}
\end{figure}

\begin{figure}[h]
\begin{center}
\includegraphics[scale=0.4]{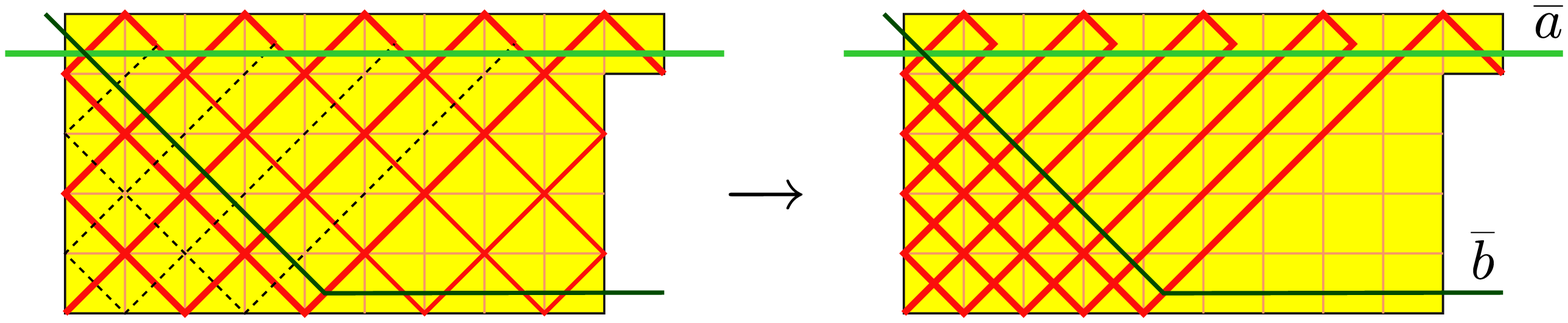}
\caption{Deformation: Case $j$ is odd \ (ex. $j = -5$)}
\label{fig:DivOdneg}
\end{center}
\end{figure}

[Case $j<-1$] \ Finally, we study the case $j < -1$. 
The method is similar.
We define an integer $j'$ by $j = -(j'+1)$ for figures.
Starting with Baker's description in Figure~\ref{fig:Baker},
we have a link $L(j)$ in Figure~\ref{fig:fromTn}.
%
%
%%%%%%%%%%%%%%%%%
%%%%%%%%%%%%%%%%%
%
%
Divide description of $L(j)$ is in Figure~\ref{fig:DivEOneg},
contrast to Figure~\ref{fig:DivEO}.
Especially, the non-trivial component 
$c(j)$ and $c(j')$ (ex. $c(-6)$ and $c(5)$) are the same curves 
presenting $T(j, j+1) = T(j'+1, j')$ but 
the other segment components are placed differently.
%
%
%%%%%%%%%%%%%%%%%
%%%%%%%%%%%%%%%%%
%
%
In Figure~\ref{fig:DivEvneg} and 
\ref{fig:DivOdneg}, we show some deformations.
Figures~\ref{fig:fromTn},
\ref{fig:DivEOneg},
\ref{fig:DivEvneg},
\ref{fig:DivOdneg}
(Case $j <-1$) 
are contrast to 
Figures~\ref{fig:fromT},
\ref{fig:DivEO},
\ref{fig:DivEv},
\ref{fig:DivOd} (Case $j>0$),
respectively.
%
%
%%%%%%%%%%%%%%%%%
%%%%%%%%%%%%%%%%%
%
%
%
%
%%%%%%%%%%%%%%%%%
%%%%%%%%%%%%%%%%%
%
%
\qed

\medskip

{\bf Proof} (of Lemma~\ref{lem:Pj}). \ 
We study a sublink
$k(j) \cup u_{-1}$ of $L(j)$
in Figure~\ref{fig:fromT} (or Figure~\ref{fig:fromTn}, respectively),
presented by $c(j) \cup \ell$ in $C(j)$ in Figure~\ref{fig:DivEO} 
(or Figure~\ref{fig:DivEOneg}) as a divide link.
The component $k(j)$ is a torus knot $T(j,j+1)$ (or $T(j'+1, j')$).
By the diagrams of $L(j)$ in Figure~\ref{fig:DivEO}
and Figure~\ref{fig:DivEOneg}, 
we can see that $\overline{k(j)}$ is a torus knot $T(j,2j+1)$
(or $T(j'+1, 2j'+1)$), presented by $P(j)$ as a divide knot.
The case $j<-1$ is a little harder, see Lemma~\ref{lem:Tn2n-1}.
We have the lemma.
\qed
\medskip

{\bf Proof} (of Lemma~\ref{lem:Pm}). \ 
In \cite{Y3}, a divide knot presentation of cable knots 
(under some conditions) is studied.
Here we use $\Delta$-moves on divides freely.

First, assume $j>0$. 
The plane curve $P_m(j)$ is obtained by adding a square twice
from $c(j)$ isotopic to $B(j,j+1)$: 
we add a square along $\ell$ to $c(j)$ first (then we have $P(j)$),
and add another square along $E_a$ or $\overline{a}$ second.
We see the plane curve obtained by the first square addition (blow-down)
in Figure~\ref{fig:DivEv} or Figure~\ref{fig:DivOd}.
Since the line $\overline{a}$ can be moved to an edge 
of the L-shaped region by $\Delta$-moves, 
we can add the second square.
The curve becomes an L-shaped curve of type $[(3j,2j),(3j+1,j)]$.
By \cite{Y3}, it presents the required cable knot $C(T(2,3); j, 6j + 1)$.

The proof in the case $j<-1$ is similar.
From the first curves in Figure~\ref{fig:DivEvneg} or Figure~\ref{fig:DivOdneg},
we have an L-shaped curve of type 
$[(\vert j \vert +1,3\vert j \vert ),(2\vert j \vert ,3\vert j \vert -1)]$.
By \cite{Y3}, we have the lemma.
\qed
\par

\bigskip
%%%%%%%%%%%%%%%%%%%%%%%%%%%%%%%%%%%%%%%%%%
{\bf Acknowledgement.}
The authors would like to thank to Professor Olivier Couture
for his valuable advice in the opportunity of a conference
^^ ^^ Singularities, knots, and mapping class groups in memory of Bernard Perron"
held in Sept.\ 2010.
Without Couture's method, the proof would be troublesome and longer.
The author would like to thank to 
Professor Mikami Hirasawa,
%Professor Tomomi Kawamura, 
%Professor Masaharu Ishikawa 
and Professor Norbert A'Campo for informing him of divide knot theory.
The author also would like to thank to 
Professor Kimihiko Motegi, 
Professor Toshio Saito, 
%Professor Motoo Tange and  
Professor Kenneth L Baker 
%and Professor John Berge
for helpful suggestions on lens space surgery.

This work was supported by JSPS KAKENHI (Grant-in-Aid for Scientific Research) (C) 
Grant Number 24540070.
%The author would like to express his
%sincere gratitude to the referee for carefully reading the manuscript and giving him
%much advice.
%%%%%%%%%%%%%%%%%%%%%%%%%%%%%%%%%%%%%%%%%%

%%%%%%%%%%%%%%%%%%%%%%%%%%%%%%%%%%%%%%%%%%
%%%                   References     Y Yamada
%%%%%%%%%%%%%%%%%%%%%%%%%%%%%%%%%%%%%%%%%%
%%%%%%%%%%%%%%%%%%%%%%%%%%%%%%%%%%%%%%%%%%
%%% END of Footnotesize

{\small
\par
YAMADA Yuichi \par
Dept. of Mathematics, 
The University of Electro-Communications \par
1-5-1,Chofugaoka, Chofu, Tokyo, 182-8585, JAPAN \par
{\tt yyyamada=AT=sugaku.e-one.uec.ac.jp} \par
}
\end{document}